\documentclass[11pt,draftcls,onecolumn]{IEEEtran}
\usepackage{amsfonts,amssymb,amsbsy,amsmath,paralist,theorem,bm,ifthen,color}
\usepackage[pdfstartview=FitH,bookmarksnumbered,unicode,bookmarksopen=true]{hyperref}
\usepackage{graphicx}
\usepackage{subfigure,relsize}
\usepackage[numbers,sort&compress,square]{natbib}
\usepackage{cases,booktabs}
\usepackage{algorithmic,algorithm}
\usepackage{array}

\newtheorem{property}{Property}

\newtheorem{rem}{Remark}

\renewcommand\sb{\ensuremath{{\bf s}}}

\newcommand\Ab{\ensuremath{{\bm A}}}
\newcommand\ab{\ensuremath{{\bm a}}}
\newcommand\Bb{\ensuremath{{\bm B}}}
\newcommand\bb{\ensuremath{{\bm b}}}
\newcommand\Cb{\ensuremath{{\bm C}}}
\newcommand\cb{\ensuremath{{\bf c}}}

\newcommand\Ib{\ensuremath{{\bm I}}}

\newcommand\pb{\ensuremath{{\bm p}}}

\newcommand\qb{\ensuremath{{\bm q}}}

\newcommand\zb{\ensuremath{{\bm z}}}

\newcommand\lambdab{\ensuremath{{\bm \lambda}}}

\newcommand\zerob{\ensuremath{{\bm 0}}}

\newcommand\diag{\ensuremath{{\rm diag}}}

\newcommand\oneb{\ensuremath{{\bf 1}}}

\newcommand\Xset  {\ensuremath{{\mathcal{X}}}}
\newcommand\Zset  {\ensuremath{{\mathcal{Z}}}}
\newcommand\Nset  {\ensuremath{{\mathcal{N}}}}
\newcommand\Mset  {\ensuremath{{\mathcal{M}}}}
\newcommand\Kset  {\ensuremath{{\mathcal{K}}}}
\newcommand\Rset  {\ensuremath{{\mathbb{R}}}}
\newcommand\st    {\ensuremath{{\rm s.t.}}}
\newcommand\Pmax  {\ensuremath{P_{\max}}}
\newcommand\dmax  {\ensuremath{d_{\max}}}
\newcommand\dmin  {\ensuremath{d_{\min}}}
\newcommand\Hmax  {\ensuremath{H_{\max}}}
\newcommand\Hmin  {\ensuremath{H_{\min}}}

\DeclareMathOperator*{\argmax}{argmax}
\DeclareMathOperator*{\argmin}{argmin}

\def\blue{\color{blue}}

\begin{document}
\title{Multi-UAV Interference Coordination via Joint Trajectory and Power Control}
\author{Chao Shen, Tsung-Hui Chang, Jie Gong, Yong Zeng, and Rui Zhang 
\thanks{\smaller[1]	C. Shen is with Beijing Jiaotong University, China (email: chaoshen@bjtu.edu.cn, chaoshen@ieee.org). }
\thanks{\smaller[1]
	T.-H. Chang is with the School of Science and Engineering, The Chinese University of Hong Kong, Shenzhen and 
	Shenzhen Research Institute of Big Data, Shenzhen, China 518172 (email: tsunghui.chang@ieee.org). }
\thanks{\smaller[1]	J. Gong is with Sun Yat-Sen University, China (email: gongj26@mail.sysu.edu.cn).}
\thanks{\smaller[1]    Y. Zeng is with the School of Electrical and Information Engineering, The University of Sydney, Australia (email: yong.zeng@sydney.edu.au).}
\thanks{\smaller[1] R. Zhang is with the National University of Singapore, Singapore (email: elezhang@nus.edu.sg).} 
}
\maketitle	

\begin{abstract}
In this paper, we consider an unmanned aerial vehicle-enabled interference channel (UAV-IC),  where each of the $K$ UAVs communicates with its associated ground terminals (GTs) at the same time and over the same spectrum.
To exploit the new degree of freedom of UAV mobility for interference coordination between the UAV-GT links,
we formulate a joint  trajectory and power control (TPC) problem for maximizing the aggregate sum rate of the UAV-IC for a given flight interval, under the practical 
constraints on the UAV flying speed, altitude, as well as collision avoidance. 
These constraints couple the TPC variables across different time slots and UAVs, leading to a challenging large-scale and non-convex optimization problem.
By exploiting the problem structure, we show that the optimal TPC solution  follows the fly-hover-fly strategy, based on which the problem can be handled by firstly finding an optimal hovering locations followed by solving a dimension-reduced TPC problem with given initial and hovering locations of UAVs. For the reduced TPC problem, we propose a successive convex approximation algorithm.
To improve the computational efficiency, we further develop a parallel TPC algorithm that is effciently implementable over multi-core CPUs. We also propose a {segment-by-segment method} 
which decomposes the TPC problem into sequential TPC subproblems each with a smaller problem dimension.
Simulation results demonstrate the superior computation time efficiency of the proposed algorithms, and also show that the UAV-IC can yield higher network sum rate than the benchmark orthogonal schemes.
\end{abstract}
\noindent {\bf Keywords} - UAV, trajectory design, power control,  interference channel, collision avoidance, distributed optimization. 

\section{Introduction}\label{sec::intro}
An unmanned aerial vehicle (UAV), commonly known as a drone, is an aerial vehicle that can fly autonomously or be piloted by a ground control station to accomplish commercial or military missions. 
In recent years, UAVs have found many promising applications in wireless communications \cite{Zeng2016Mag,FlyingBS} due to their high agility, ability of on-demand and low-altitude deployment, and strong communication links with the ground.
For example, UAVs can be deployed rapidly as aerial base stations or aerial mobile relays to provide enhanced communication performance for existing wireless communication 
networks or support emergent service in war or disaster areas. Besides, UAVs can be used to perform remote surveillance and deliver real-time video data to ground terminals (GTs) \cite{Motlagh2017}. 
UAVs are also useful for data collection and dissemination in wireless sensor networks  \cite{T:MAC:1UAVzhan,UAV4Dissemination,Jie2018}.

Unlike the conventional wireless systems with static access points, UAV-enabled wireless communications often require joint consideration of trajectory planning and communication resource allocation. This is because UAVs are usually energy constrained, hence optimal trajectory planning for saving energy is of great importance in  UAV applications \cite{T:P2Pzeng}. 
Besides, interference management is challenging in UAV-enabled multi-user wireless communications.  Since the air-to-ground (A2G) channel between UAV and GT usually consists of a strong line-of-sight (LoS) link \cite{3GPP_A2G}, strong interference may be inevitable for both uplink and downlink transmissions, which is in a sharp contrast to the conventional terrestrial wireless networks with severe shadowing and channel fading. 
While interference can be avoided by orthogonal transmission schemes, e.g., time-division multiple-access (TDMA) and frequency-division 
multiple-access (FDMA), they usually lead to suboptimal and low spectrum utilization efficiency.
Therefore, it is desirable to develop new interference management techniques that cater to the unique A2G channel characteristic and leverage the controllable mobility of UAVs. 
However, UAVs are usually subject to many flying constraints, such as no limited zone, limited flight speed, minimum and maximum flying altitude and so on.
When there are multiple UAVs,  the maintenance of safe separation between UAVs, i.e., collision avoidance (CA), is another critical challenge \cite{TCST2007} in the UAV trajectory design. 

\subsection{Related Work}
There have been growing research efforts for UAV-enabled wireless communication systems.
References \cite{OD2013,OD2015,OD201708,OD2017Lyu,OD2017Alzenad,JSAC_dep} considered the optimal placement of UAVs for either providing guaranteed quality-of-service (QoS) for fixed GTs or maximizing the service coverage over a given area. 
Reference \cite{T:Relay:1S1R1Dzeng} studied a UAV-enabled mobile relaying system, and proposed an efficient algorithm for joint UAV trajectory and power control (TPC). 
The authors in \cite{T:P2Pzeng} studied the energy-efficient UAV trajectory optimization problem by taking into account the UAV's propulsion energy consumption.
Reference \cite{T:MAC:1UAVzhan} considered the UAV trajectory design for data collection in wireless sensor networks, 
and the authors in \cite{T:BC:1UAVlyu} studied an interesting throughput-delay trade-off for UAV-enabled multi-user system.
Reference \cite{T:MAC:1UAVswindlehurst} considered the UAV heading optimization for an uplink scenario with multiple antennas at the UAV. 
A general UAV-enabled radio access network (RAN) supporting multi-mode communications of the ground users was considered in \cite{RZhangTSP2018}, where new designs for the UAV initial trajectory were proposed.
The work \cite{T:Relay:KS1RKDsineone} utilized the UAV to offer dynamic computation offloading for GTs. Besides, a UAV-enabled wireless power transfer (WPT) system is studied in {\cite{Xu_UAV4WPT}}, where the harvested energy profile of GTs is characterized by optimizing the UAV trajectory.

Besides the works considering single UAV only, references \cite{T:MAC:NUAVgunduz,Wu2018TWC, UAVCoMP} studied more complicated scenarios with multiple UAVs.
Specifically, \cite{T:MAC:NUAVgunduz} considered a scenario with multiple UAVs communicating with one GT and derived the optimal TPC for minimizing both transmission and propulsion energy.
 \cite{Wu2018TWC} considered the joint TPC and user association/scheduling problem for multiple UAVs serving multiple GTs. 
 \cite{UAVCoMP} considered the use of multiple UAVs for coordinated multipoint (CoMP) transmission to serve a set of GTs, and derived the optimal TPC and UAV deployment solutions for maximizing the ergodic sum rate of users under random channel phases.

\subsection{Contributions}

In this paper, we study a UAV-enabled interference channel (UAV-IC), where $K$ UAVs communicate with their respective GTs at the same time and over the same spectrum,  as shown in Fig. \ref{fig:fig1}. 
This scenario is well motivated in practice. For example, one use case is that the UAVs collect data from a field and deliver the data to their respective serving GTs \cite{T:MAC:1UAVzhan,UAV4Dissemination}. Note that 
the UAV-IC considered in this paper is generic and the techniques developed here can be generalized to other similar scenarios with multiple UAVs \cite{T:MAC:NUAVgunduz,Wu2018TWC, UAVCoMP}.

In contrast to traditional interference channels with static terminals \cite{NPhardness,WMMSE}, the UAV-IC allows one to exploit the mobility of UAVs as a new degree of freedom to dynamically control the interference among $K$ communication links via joint trajectory optimization and power control. 
Therefore, we formulate a joint TPC problem for maximizing the aggregate sum rate of all UAV-GT pairs. Different from most prior works on UAV trajectory optimization, we consider not only the practical constraints on the flying speed and altitude of each UAV, but also the minimum spacing constraint between UAVs for collision avoidance in the three-dimensional (3D) space. 

The formulated TPC problem for the UAV-IC has several major challenges. 
Firstly, the problem is NP-hard in general since the sum rate maximization problem in interference channels even with static terminals is NP-hard  \cite{NPhardness}. 
Secondly, the formulated problem involves joint TPC for all the $K$ UAVs and for each UAV, it usually includes a large number of TPC variables due to a long flying interval.  In the existing works such as  \cite{T:Relay:1S1R1Dzeng,T:MAC:1UAVzhan, T:P2Pzeng,Wu2018TWC}, the TPC optimization is usually handled by employing the alternating optimization (AO) technique, which optimizes the trajectory  variables and the transmission power variables separately and iteratively.
Both the subproblems for trajectory optimization and power control are non-convex, and are often handled by the successive convex approximation (SCA) technique  
\cite{SCA_1,SCA_2,WMMSE}. However, the AO based methods may converge to a non-desirable local point, especially when the variables are coupled with each other in the constraints \cite{BSUM}; they may not be time efficient either since the algorithm involves two loops of optimizatation procedure (outer AO loop and inner SCA loop).

In this paper, we propose new and computationally efficient algorithms to solve the TPC problem for the UAV-IC. The main contributions are summarized below.
\begin{enumerate}
	\item  Assuming that each UAV has to return to its initial location at the end of flight, we show that the optimal solution to the TPC optimization problem has an interesting symmetric property.
	By this property, the TPC problem can be decomposed to firstly finding the optimal hovering locations of UAVs that maximize the instantaneous sum rate (i.e., the deployment problem), followed by optimizing the TPC from the initial locations to the optimal hovering locations. Both problems have significantly reduced problem dimensions than the original TPC problem.
	Similar decomposition strategies are applicable to scenarios where the initial and destination locations of UAVs are different.
	
	\item For the TPC optimization, we propose a new SCA algorithm with a locally tight surrogate rate function that allows joint update of the trajectory variables and transmission powers in each iteration, which is different from the AO based algorithms in \cite{T:Relay:1S1R1Dzeng,T:MAC:1UAVzhan, T:P2Pzeng,Wu2018TWC}.
	Numerical results show that the proposed algorithm requires less computation time than the AO based method, while both of them can achieve compariable sum rate performances.
	
	\item To overcome the high computational complexity due to the large number of UAVs, based on the block successive upper bound minimization method of multipliers (BSUMM) \cite{BSUM}, we propose a distributed TPC algorithm that can be implemented in parallel over multi-core CPUs. 
	To further reduce the computation time when the flight duration is large, we propose a segment-by-segment strategy, which divides the entire flight trajectory 
	into consecutive smaller time segments. 
	The TPC optimization for each time segment is thus more efficiently solvable.
\end{enumerate}
The proposed SCA technique is also applied to the TPC problems when the UAVs adopt FDMA or TDMA schemes. 
Simulation results not only show the superior computational efficiency of the proposed algorithms over the conventional AO methods and improved spectrum efficiency than FDMA/TDMA. 

The rest of this paper is organized as follows. Section \ref{sec: system model} presents the system model of the UAV-IC and the corresponding TPC optimization problem for maximizing the aggregate sum rate.
The TPC optimal solution is analyzed and 
a new SCA-based algorithm for efficient TPC optimization is presented in Section \ref{sec: TPC solution}.
In Section \ref{sec: distr opt}, a distributed TPC algorithm and a segment-by-segment method 
are proposed for computation time reduction.
In Section \ref{sec: FDMA/TDMA}, the proposed SCA algorithms are extended to FDMA and TDMA schemes.
The simulation results are given in Section \ref{sec:simulation} and the work is concluded in Section \ref{sec: conclusion}.

{\bf Notations:} Column vectors and matrices are respectively written in boldfaced lower-case and upper-case letters, e.g., $\ab$ and $\Ab$. The superscript $(\cdot)^T$ represents the transpose.
$\Ab\succeq \zerob$ means that matrix $\Ab$ is positive semidefinite. $\Ib_K$ is the $K \times K$ identity matrix; $\oneb_K$ is the $K$-dimensional all-one vector;  $\zerob_K$ is the $K$-dimensional all-zero vector. $\|\ab\|_2$ denotes the Euclidean norm of vector $\ab$, and $\|\zb\|^2_\Ab\triangleq \zb^T\Ab\zb$ for some $\Ab\succeq \zerob$.
$\diag(\{a_i\})$ and $\diag(\ab)$, where $\ab=[a_1,\ldots,a_m]^T$, both represent a diagonal matrix with $a_i$'s being the diagonal elements.
Notation $\otimes$ denotes the Kronecker product.


\section{System Model and Problem Formulation}\label{sec: system model}
\begin{figure}\centering
	\includegraphics[width=0.7\linewidth]{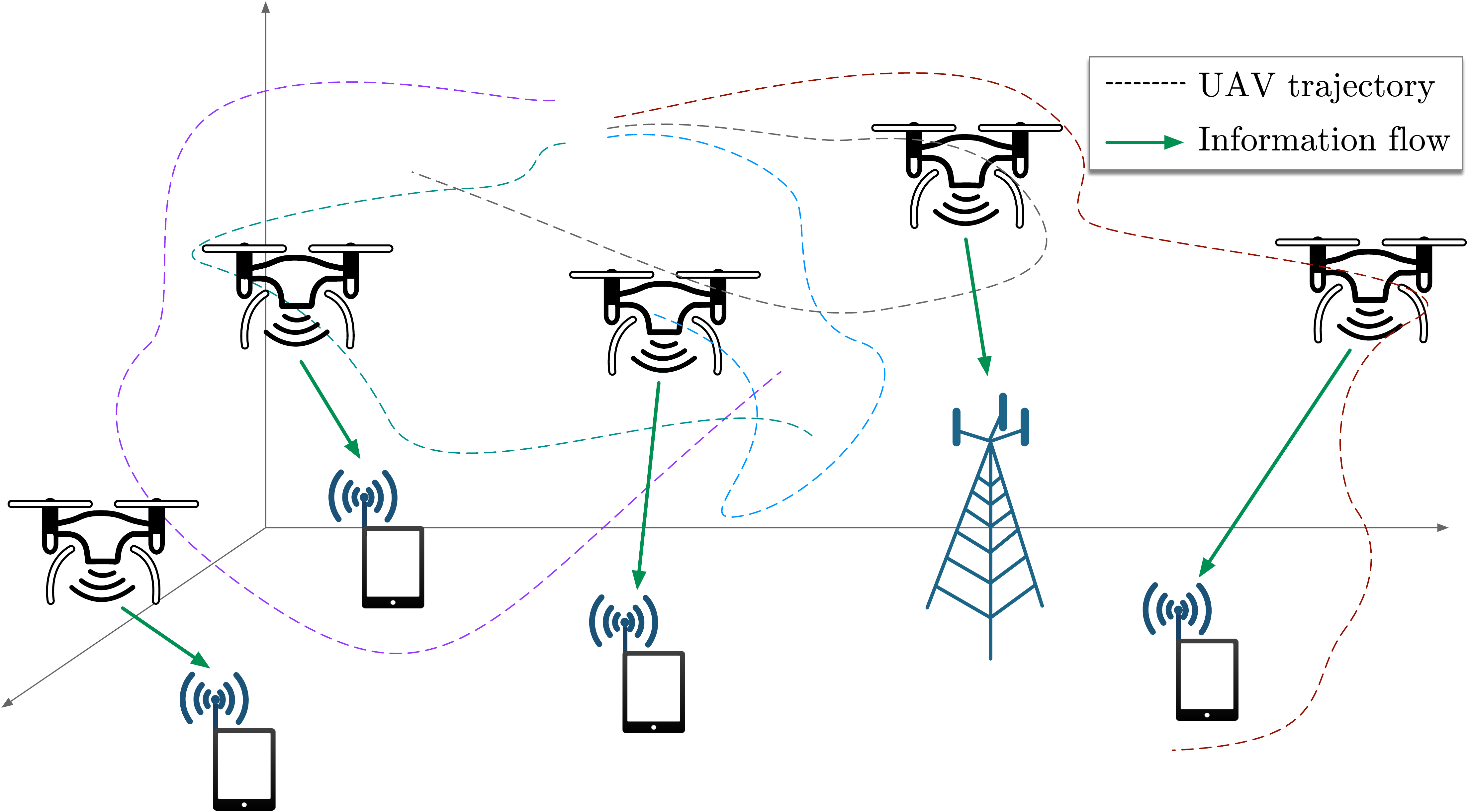}
	\caption{ The considered UAV-IC network with $K=5$ UAV-GT pairs in a 3D Cartesian coordinate system. 
		The UAVs collect data from the field and deliver the data to their respective serving GTs.}
	\label{fig:fig1}
\end{figure}

As shown in Fig. \ref{fig:fig1}, we consider a UAV-enabled wireless communication system where $K$ UAVs respectively communicate with their associated GTs at the same time and over the same spectrum.  For convenience, we consider the downlink communication from the UAVs to GTs only, though the developed results also apply to the uplink communication as well.  
Both UAVs and GTs are assumed to be equipped with one omnidirectional antenna.
The GT locations, denoted by $\sb_k\in\Rset^3$ for $k\in\Kset\triangleq \{1,\ldots,K\}$, are assumed to be known.
Our goal is to jointly optimize the flight trajectories of the $K$ UAVs and the transmission powers for a given time horizon $T$ 
so as to maximize the network throughput. 

For ease of design, the time horizon $T$ is discretized into $N$ equally spaced time slots, i.e., $T=NT_s$, with $T_s$ being the sampling interval.
For each $k\in\Kset$, denote $\qb_k[n]\in\Rset^3$ as the location of the $k$th UAV at time slot $n$, for $n\in \Nset\triangleq \{1,\ldots,N\}$. We denote $\qb_k[0]$ and $\qb_k[N+1]$ as the initial and final locations of the $k$th UAV, respectively. 

Moreover, due to the mechanical and regulatory limitations, $\{\qb_k[n]\}$ are subject to several constraints.  
Firstly, a UAV is usually subject to the minimum and maximum altitude constraints, which can be represented as\vspace{-2pt}
\begin{align}\label{altitute constraint}
\Hmin \leq (\qb_k[n])_3 \leq \Hmax,  
\end{align}
for all $n\in\Nset, ~k\in\Kset$, where $(\qb_k[n])_3$ represents the third element of $\qb_k[n]$. 
Secondly, the flight speed of each UAV is typically limited by the maximum level-flight speed $V_L$, vertical ascending speed $V_{A}$  and vertical descending speed $V_{D}$. 
For example, $V_L=17$ m/s,  $V_{A}=5$ m/s, and  $V_{D}=3$ m/s for DJI M200 \cite{DJIm200}, which is a commercial quadrotor used for power line inspection. 
The maximum UAV speed constraints are equivalent to constraining the maximum distances that a UAV can move during one time interval, i.e.,\vspace{-3pt}
\begin{subequations}\label{FullSpeedLimit}
	\begin{align}
	\!\!\!\|(\qb_k[n])_{1:2}-(\qb_k[n-1])_{1:2}\| & \leq d_L\!\triangleq\! V_L T_s,\\
	\!\!\!-V_{D} T_s\!\triangleq\! -d_D \leq (\qb_k[n])_3\!-\!(\qb_k[n-1])_3  &\leq d_A\!\triangleq\! V_{A}T_s,
	\end{align}
\end{subequations}
for all $k\in\Kset$, and $n\in\bar\Nset\triangleq\Nset\cup\{N+1\}$, where $(\qb_k[n])_{1:2}$ denotes the first two components of the vector $\qb_k[n]$.
Thirdly, to ensure collision avoidance, any two UAVs need to be separated by a minimum distance $\dmin$ at any time instance, that is, for all $j$, $k\in \Kset$, $j>k$ and $n\in \Nset$, 
\begin{align}\label{eqn:CA constraint}
\|\qb_j[n] - \qb_k[n]\| \geq\dmin.
\end{align} 
 

Previous channel field measurements have shown that the LoS component dominates the A2G channels in many practical scenarios \cite{3GPP_A2G}, especially for rural areas or moderately high UAV altitude. Therefore, in this paper,  the channel gain between $j$th UAV and $k$th GT 
is modeled by the free-space path loss model, which can be expressed as $\frac{\beta_0}{\|\qb_j[n] - \sb_k\|^2}$.  Here, $\beta_0$ denotes the channel power gain at the reference distance of one meter. Thus, for the considered $K$-user UAV-IC, the achievable rate in bits/sec (bps) of the $k$th UAV-GT pair at time slot $n$ is given by 
\begin{align}\label{Rate}\!\!
R_k(\pb[n],\qb[n])\!\triangleq\! B\log_2\!\!\Bigg(\!1\!+\! \frac{\frac{\gamma p_k[n]}{\|\qb_k[n]-\sb_k\|^2}}{ 1+\sum_{j=1,j\neq k}^K\frac{\gamma p_j[n]}{\|\qb_j[n]-\sb_k\|^2}}\!\Bigg),
\end{align}
where $\pb[n]\!\triangleq\! [p_1[n],\!\ldots,\!p_K[n]]^T$ and $\qb[n]\!\triangleq\![\qb_1^T[n],\!\ldots,\!\qb^T_K[n]]^T$ contain the transmission powers and locations of all UAVs at time slot $n$, and $\gamma\triangleq \frac{\beta_0}{BN_0}$ with $N_0$ denoting the power spectral density (PSD) of the additive white Gaussian noise (AWGN) and $B$ being the communication bandwidth. As seen from \eqref{Rate}, the interference between the UAV-GT links depends  on both the UAV trajectories and transmission power. 

 Our goal in this paper is to jointly optimize the trajectories and transmission powers of all the UAVs so as to maximize the 
the aggregate sum rate of all the UAV-GT pairs over the entire flight duration. 
Such a joint TPC design problem can be formulated as follows
\begin{subequations}\label{Problem::SEM}
	\begin{align}	\!\!
	{\sf (TPC)}\max_{\substack {p_k[n],\qb_k[n]\\ k\in\Kset,n\in\Nset} }&~\sum_{n=1}^N \sum_{k=1}^{K}  R_k(\pb[n],\qb[n])  \\
	\st~~~	
	&\!\!\!\!\eqref{altitute constraint}, \eqref{FullSpeedLimit}, \eqref{eqn:CA constraint}~ \forall n\in\bar\Nset{\rm~and~} j,k\in \Kset, j\!>k, \label{eqn::SEM C0}\\
	&0 \leq p_k[n]\leq P_{\max},				~\forall k\in \Kset, n\in\Nset,  					\label{eqn::SEM C1}\vspace{-.5\baselineskip}
		\end{align}
\end{subequations}
where $P_{\max}$ is the maximum transmission power of the UAVs. 
In \eqref{Problem::SEM}, the initial and final UAV locations of UAV $k$ are fixed at $\qb_k[0]$ and $\qb_k[N+1]$, respectively, for all $k$.
Specifically, we assume that $\qb_k[0]=\qb_k[N+1]$ for all $k$, i.e., each UAV has to return to its initial location at the end of the flight. 
This arises, for example, in wireless sensor networks where the UAVs 
are collecting sensory data around the locations $\{\qb_k[0]\}$ and have to deliver the data to their respective GTs followed by returning back to $\{\qb_k[0]\}$ for further data collection \cite{Zeng2016Mag,FlyingBS,Kobayashi2014,Motlagh2017,T:MAC:1UAVzhan,UAV4Dissemination}.

We should emphasize that the TPC problem \eqref{Problem::SEM} is challenging to solve. Firstly, even for the conventional $K$-user interference channel with static terminals, the power allocation problem for sum rate maximization is known to be NP-hard for $K>2$ \cite{NPhardness}. Therefore, the TPC problem \eqref{Problem::SEM} is also NP-hard.
Secondly, the problem dimension of \eqref{Problem::SEM} can be large. Notice that the sampling interval $T_s$ should be small enough (i.e.,  the number of discrete time slots $N$ should be large enough) so that 
the channel can be treated as approximately constant within each time slot, and furthermore, the collision between UAVs can be effectively avoided. Specifically, it requires 
\begin{align}\label{maxTs}
T_s \leq \frac{\dmin}{\sqrt{4V_L^2 + (V_D + V_A)^2}}
\end{align}
so that in the worst case that two UAVs fly toward each other with the maximum speed, collision between the two UAVs can be effectively detected and avoided.
In light of these computational challenges, it is desirable to develop computationally efficient algorithms for obtaining high-quality suboptimal solutions to the TPC problem \eqref{Problem::SEM}.


\section{Proposed TPC Optimization Algorithms}\label{sec: TPC solution}
In Section \ref{subsec: TPC solution property}, we show that the TPC problem \eqref{Problem::SEM} can be handled by first finding an optimal deployment location for each UAV, followed by optimizing the trajectories and transmission powers of the UAVs from their initial locations to the deployment locations. 
%
This enables us to consider a dimension-reduced counterpart of problem \eqref{Problem::SEM} as shown in Section \ref{subsec: 2-stage approx}.
Then in Section \ref{subsec: SCA}, a novel SCA-based algorithm for the TPC design problem is presented.\vspace{-0.7\baselineskip}

\subsection{Symmetry of TPC Solution}\label{subsec: TPC solution property}
The optimal solution to the TPC problem \eqref{Problem::SEM} is not unique in general. 
Interestingly, under the round-trip constraint of $\qb_k[0]=\qb_k[N+1]$ for all $k\in\Kset$, the TPC solution is symmetric with respect to the flight trajectories.
Without loss of generality, we assume that $N$ is even.\vspace{-0.6\baselineskip}
\begin{property}\label{Prop:Rmax}
	The TPC problem \eqref{Problem::SEM} admits an optimal solution satisfying
	\begin{align}\label{Prop:eq:symmetric}
	   \qb_k^\star[n]=\qb_k^\star[N+1-n], ~p_k^\star[n]=p_k^\star[N+1-n],
	\end{align}
	for all $n\in  \Nset$, $k\in \Kset$. 
Moreover, for some $M\in\{1,2,\ldots,N/2\}$,
	\begin{align}
	\!\!R_s[n] \!\leq\! R_s[M]\!=\!R_s[M+1]\!=\!\cdots \!=\!R_s[N/2],~\forall n\!\in\!\Mset. \label{Prop:eq:Rmax}
	\end{align}
	where $R_s[n] \triangleq \sum_{k=1}^K R_k(\pb^\star[n],\qb^\star[n])$ is the sum rate at time slot $n$, and $\Mset\triangleq \{1,\ldots,M\}$. \vspace{-0.6\baselineskip}
\end{property}
\begin{IEEEproof}
	Let $R_1=\sum_{n=1}^{N/2} R_s[n]$ and $R_2=\sum_{n={N/2}+1}^{N} R_s[n]$ be the aggregate sum rates in the first half and second half of the time horizon, respectively.
	Suppose that \eqref{Prop:eq:symmetric} is not true and without loss of generality (w.l.o.g.), $R_1>R_2$.
	Then simply assigning the trajectory and transmission power solutions of UAVs at time slot $n$ in the first half flight to those at time slot $N+1-n$ in the second half flight for all $n\in [1,N/2]$ can achieve an objective value of $2R_1$, which is strictly larger than $R_1+R_2$.  Furthermore, it is not difficult to show that this new trajectory and power allocation satisfy all constraints in {\sf (TPC)}. This contradicts with the optimality of  $\{\pb^\star[n],\qb^\star[n]\}_{n=1}^N$. Thus, the optimal TPC solution must have $R_1=R_2$ which admits 	
	a symmetric TPC solution as in \eqref{Prop:eq:symmetric}.
	
	We then focus on the TPC solutions for the first half flight.
	Suppose that \eqref{Prop:eq:Rmax} is not true.
	Due to the constraints in \eqref{eqn::SEM C0} and \eqref{eqn::SEM C1}, $R_s[n]$ is bounded for all $n\in{\Nset}$, and thus there exists  some $M\in\{1,2,\ldots,N/2\}$ such that $R_s[M]=\max \{R_s[1],\ldots,R_s[N/2]\}$. Then, we must have
		\begin{align}\label{Vmodel:eq2}
	\left(N/2-M+1\right)R_s[M]\ge \sum_{n = M}^{N/2} R_s[n],
	\end{align}
	which implies that by assigning the trajectory and power solutions of UAVs at time slot $M$ to those in time slots $M+1,\ldots,N/2$, we can achieve an objective value of $\sum_{n = 1}^{M-1} R_s[n] + \left(N/2-M+1\right)R_s[M] \ge \sum_{n = 1}^{N/2} R_s[n]$ without violating any constraints in {\sffamily(TPC)}.
    This is a contradiction and shows that \eqref{Prop:eq:Rmax} is true.
\end{IEEEproof}

Property \ref{Prop:Rmax} is insightful. 
It implies that the \emph{fly-hover-fly} strategy is optimal, that is, the UAVs should first fly to the locations $\{\qb_k[M]\}$ which have the maximum sum rate $R_s[M]$ along the trajectories, hover over the locations until time slot ${N}-M$, and then fly back to the initial locations along the same paths. Hence, the TPC problem \eqref{Problem::SEM} can be reduced to the TPC design problem for the first $M$ time slots only, rather than for the entire duration with $N$ time slots.
Specifically, it is equivalent to solve the following dimension-reduced counterpart of \eqref{Problem::SEM}\vspace{-0.4\baselineskip}
\begin{subequations}\label{Problem::SEM2}
	\begin{align}	
	&\max_{\substack {p_k[n],\qb_k[n]\\ k\in\Kset,n\in \Mset \\ M \in \{1,2,\ldots,N/2\} } }\sum_{n=1}^M \sum_{k=1}^{K}  R_k(\pb[n],\qb[n])\notag\\[-18pt]
	&~~~~~~~~~~~~~~~~~~~~~~+ \left(\tfrac{N}{2}-M\right)  \sum_{k=1}^{K}  R_k(\pb[M],\qb[M])   	\label{eqn::SEM2 obj2} \\
	&~~~~~~~\st~~~ \eqref{altitute constraint}, \eqref{FullSpeedLimit}, \eqref{eqn:CA constraint}~\forall n\in\Mset{\rm~and~} j,k\in \Kset, j>k, \label{eqn::SEM2 C0}\\
&~~~~~~~~~~~~~~ 0 \leq p_k[n]\leq P_{\max},~\forall k\in \Kset, n\in\Mset.  					\label{eqn::SEM2 C1}\vspace{-1\baselineskip}
	\end{align}
\end{subequations}

\subsection{Heuristic 2-Step TPC Design}\label{subsec: 2-stage approx}

Firstly, note that the optimal solution structure in \eqref{Prop:eq:Rmax} is still applicable to problem \eqref{Problem::SEM2}.
Secondly,  in \eqref{Problem::SEM2}, we need to optimize the value of $M$, which however is difficult in general. 
In this paper, we adopt a heuristic 2-stage procedure to handle problem \eqref{Problem::SEM2}: the first stage determines the hovering locations $\{\qb_k[M]\}$ and the value of $M$, and the second stage solves \eqref{Problem::SEM2} with fixed $\{\qb_k[M]\}$ and $M$.
The procedure is motivated by the observation that the second term of \eqref{eqn::SEM2 obj2} would become dominant if $N$ is much larger than $M$\footnote{By intuition, if $N$ is much larger than the number of time slots required for the UAVs flying to the GTs, we would have $N\gg M$.}.
In that case, the hovering locations $\{\qb_k[M]\}$ would be approximately the same as the locations 
where the UAVs achieve the maximum instantaneous sum rate.
Thus, we determine the hovering locations $\{\qb_k[M]\}$
by solving the following deployment problem\vspace{-0.3\baselineskip}
\begin{subequations}\label{Problem::DP}
	\begin{align}	
\!\!\!\!\!\!\!\!\!\!(\pb^\star,\qb^\star) =& \argmax_{\pb,\,\qb}~ \sum_{k=1}^{K}  R_k(p_k,\qb_k)  \\
	\st~
	&\Hmin \leq (\qb_k)_3 \leq \Hmax,		~k\in \Kset,					\label{Problem::DP C4}\\
	&\|(\qb_k)_{1:2}-(\qb_k[0])_{1:2}\| \leq V_L T/2,~k\in \Kset,            \label{Problem::DP C2-H}
	\end{align}
	\begin{align}
	&-V_D T/2 \leq (\qb_k)_3-(\qb_k[0])_3  \leq V_A T/2,~k\in \Kset,				\label{Problem::DP C2-V}\\
	&\|\qb_k-\qb_j\|  \geq d_{\min}, 	~k,j \in \Kset, j>k,  				\label{Problem::DP C3}\\
	&0 \leq p_k\leq P_{\max},				~k\in \Kset, 					\label{Problem::DP C1}
	\end{align}
\end{subequations}
where $\pb\!\triangleq\! [p_1,\!\ldots,p_K]^T$, $\qb\!\triangleq\! [\qb_1^T,\!\ldots,\!\qb_K^T]^T$, and the constraints \eqref{Problem::DP C2-H}, \eqref{Problem::DP C2-V} guarantee that the hovering positions can be reached by the UAVs within $N/2$ time slots. Different from problem \eqref{Problem::SEM2}, the deployment problem \eqref{Problem::DP} does not involve the trajectory design of the UAVs but simply determines the optimal location that achieves the maximum instantaneous sum rate, among all feasible locations where the UAVs are reachable from their initial locations within the  duration $T/2$. 

Let $\qb_k[M]=\qb_k^\star$ and $p_k[M]=p_k^\star$  for all $k\in \Kset$ be fixed in \eqref{Problem::SEM2}, and denote $R_s^\star = \sum_{k=1}^K R_k(\pb^\star,\qb^\star)$. 
We arrive at the following TPC design problem
	\begin{subequations}\label{Problem::SEM3}
	\begin{align}	
	\{\qb^\star[n],\pb^\star[n]\}_{n=1}^M =&\argmax_{\substack {p_k[n],\qb_k[n]\\ k\in\Kset,n\in \Mset} }~\sum_{n=1}^M \sum_{k=1}^{K}  R_k(\pb[n],\qb[n]) +(\tfrac{N}{2} -M )R_s^\star \\
	&~~~~~~~~~\st~~~
	 \eqref{altitute constraint}, \eqref{FullSpeedLimit}, \eqref{eqn:CA constraint} {~\rm for~}n\in  \Mset, k,j\in \Kset,j>k, \\
	&~~~~~~~~~~~~~~~~ 0 \leq p_k[n]\leq P_{\max},				~k\in \Kset, n\in \Mset,					\label{eqn::SEM3 C12}\\
	&~~~~~~~~~~~~~~~~\qb_k[M]=\qb_k^\star, ~p_k[M]=p_k^\star,~\forall k\in \Kset. \label{eqn::SEM3 C42}
	\end{align}
\end{subequations}

\begin{property}\label{Prop:M}
	The optimal objective value of \eqref{Problem::SEM3}  is non-increasing in $M$.
\end{property}
\begin{IEEEproof} Suppose that \eqref{Problem::SEM3} is feasible for some values of $M$ and $M+1$.
	Let $R_s[n] = \sum_{k=1}^K R_k(\pb^\star[n],\qb^\star[n])$ for $n\in \Mset$ and $R_s^\star = R_s[M]$. The corresponding aggregate sum rate is given by
	$R^\star \triangleq  \sum_{n=1}^M R_s[n] + (N/2-M)R_s^\star$.
	Let $\{\pb'[n],\qb'[n]\}_{n=1}^{M+1}$ be an optimal solution of \eqref{Problem::SEM3} when $M$ is replaced by $M+1$, and $R_s'[n] = \sum_{k=1}^K R_k(\pb'[n],\qb'[n])$ for $n\in  \Mset\cup \{M+1\}$.
	Due to \eqref{eqn::SEM3 C42}, $R_s'[M+1] = R_s^\star$.
	Thus, the corresponding aggregate sum rate  is given by 
	\begin{align}
	R' &\triangleq  \sum_{n=1}^M R_s'[n] + R_s'[M+1] +(N/2-M-1)R_s^\star \notag \\
	&= \sum_{n=1}^M R_s'[n] + (N/2-M)R_s^\star \leq R^\star.
	\end{align}
	The property is proved.   
\end{IEEEproof}

Property \ref{Prop:M} implies that, once the hovering locations $\qb^\star$ are given, the UAVs should fly to the hovering locations $\qb^\star$ as quickly as possible.
 For example,  one may set
\begin{align}\label{eqn: M}
 M=\max_{k\in \Kset}  \bigg\lceil\frac{\|(\qb_k^\star)_{1:2} - (\qb_k[0])_{1:2}\|}{V_LT_s}\bigg\rceil +\triangle
\end{align}
for some $\triangle\geq 0$. In the right hand side (RHS) of \eqref{eqn: M}, the first term is the required number of time slots for the UAVs flying straight to the hovering locations with the maximum speed, and the second term $\triangle$ stands for some additional number of time slots for guaranteeing no collision between UAVs.
In practice, one may 
estimate $M$ through building feasible trajectories of UAVs from the initial locations to 
the hovering locations $\qb^\star$; details will be given in Sec. \ref{Sec::TrajInit}.

In summary, we have the following procedures to obtain a TPC solution to problem  \eqref{Problem::SEM}:
\begin{itemize}\it  
	\item [{1)}] Find an optimal hovering solution $(\pb^\star,\qb^\star)$ by \eqref{Problem::DP}.
	\item [{2)}] Solve the TPC problem \eqref{Problem::SEM3} with a given value of $M$. 
	\item [{3)}] Construct the TPC solution for the entire flight by \eqref{Prop:eq:symmetric}.
\end{itemize}

The above analysis enables us to focus on developing efficient algorithms for the TPC design problem \eqref{Problem::SEM3}. 
Since the deployment problem \eqref{Problem::DP} is a special instance of \eqref{Problem::SEM3} (with $M=1$), any algorithm developed for the latter problem can also be applied to the former.\vspace{-0.6\baselineskip}

\begin{rem}\rm
	If the UAV does not have to return to its initial location, i.e., $\qb_k[N+1]\neq \qb_k[0]$, the TPC solutions of \eqref{Problem::SEM} are not symmetric in general.
	Nevertheless, the fly-hover-fly strategy is still applicable.
	To find proper hovering locations of UAVs, one may consider the following problem
	\begin{subequations}	\label{OnewayProblem::DP}
		\begin{align}
		&\!\!\!\!\!\!\!(\pb^*,\qb^*,\tau^*) \in
		\argmax_{\pb,\,\qb,\,\tau\in[0,T]}~ \sum_{k=1}^{K}  R_k(p_k,\qb_k)  \\
		\st
		&~~{\rm constraints~ in~ } \eqref{Problem::DP C4}, \eqref{Problem::DP C3}, \eqref{Problem::DP C1},    \notag \\
		&~~ \|(\qb_k)_{1:2}\!-\!(\qb_k[0])_{1:2}\| \!\leq\!  \tau V_L ,~k\in \Kset,  \notag \\
		&~~ -\!\tau V_D  \leq (\qb_k)_3-(\qb_k[0])_3  \!\leq\! \tau V_A,~k\in \Kset,    \notag\\
		&~~ \|(\qb_k)_{1:2}\!-\!(\qb_k[N+1])_{1:2}\| \!\leq\! (T-\tau) V_L,~k\in \Kset,    \label{OnewayProblem::DP C2-H2}\\
		&~~ -\!(T\!\!-\!\tau) V_D \!\leq\! (\qb_k)_3\!-\!(\qb_k[N\!+\!1])_3  \!\leq\! (T\!-\!\tau) V_A,k\!\in\! \Kset,   \label{OnewayProblem::DP C2-V2}\vspace{-0.8\baselineskip}
		\end{align}
	\end{subequations}
	where $\tau\in[0,T]$ is the time for the UAVs flying to the hovering locations; constraints in \eqref{OnewayProblem::DP C2-H2}-\eqref{OnewayProblem::DP C2-V2} ensure that the UAVs cam reach their final locations in time $T$.
    As long as $T$ is much larger than the time required for the UAVs flying from the initial locations to the hovering locations plus the time flying from the hovering locations to the final locations, the UAVs would spend most of the time hovering and the fly-hover-fly strategy would be approximately optimal.	\vspace{-0.5\baselineskip}
%
\end{rem}

\vspace{-0.3cm}
\subsection{SCA-based TPC Optimization}\label{subsec: SCA}
The TPC design problem in \eqref{Problem::SEM3} involves optimization of both trajectory variables $\{\qb[n]\}$ and transmission powers $\{\pb[n]\}$.
Following the idea of the SCA method \cite{SCA_1,SCA_2}, we successively solve a convex approximation counterpart of the non-convex problem \eqref{Problem::SEM3}. 
In particular, suppose that at the $r$th iteration, the location and transmission power of the $k$th UAV are denoted by $\qb_k^r[n]$ and $p_k^r[n]$, respectively, for $k\in \Kset$, $n\in \Mset$. 
Denote 
\begin{align}
    a^r_k[n]&= \sqrt{p^r_k[n]}, ~d^r_{jk}[n] = \|\qb^r_j[n]-\sb_k\|^2, \\
	a_k[n] &= \sqrt{p_k[n]},~d_{jk}[n] = \|\qb_j[n]-\sb_k\|^2. \label{Eq11}
\end{align}
The interference power imposed at the $k$th GT at the $n$th time slot can be written as
\begin{align}\label{Eq12}
I_k[n]&
\triangleq \sum_{j=1,j\neq k}^K \frac{\gamma p_j[n]}{\| \qb_j[n]-\sb_k\|^2}
=\sum_{j=1,j\neq k}^K \frac{\gamma a_j^2[n]}{d_{jk}[n]}.
\end{align}
By slightly abusing the notation of $R_k$, the achievable rate of the $k$th UAV-GT pair in \eqref{Rate} can be expressed as a function in $(\ab[n],\qb[n])$, as
\begin{align} \label{ratefunction}\!\!\!
&R_k(\ab[n],\qb[n]) \notag \\
&\!=\!\log\bigg(\!1\!+\!  \sum_{j=1}^K\!     \frac{\gamma a_j^2[n]}{d_{jk}[n]} \!\bigg)
\!-\!\log\bigg(\!1\!+\!  \sum_{j\neq k}^K\! \frac{\gamma a_j^2[n]}{d_{jk}[n]} \!\bigg),
\end{align}
in nats/s/Hz, where $\ab[n]=[a_1[n],\ldots,a_K[n]]^T\in\Rset^K$. 

Note that $\frac{x^2}{y}$ is convex for $x\in \mathbb{R}, y>0$, $-\log(1+x)$ is convex for  $x>-1$, and $x^2$ is convex for all $x\in \mathbb{R}$.
So, they respectively have global linear lower bounds as
$\frac{x^2}{y}\geq 
\frac{2\bar x}{\bar y}x -\frac{\bar x^2}{\bar y^2} y$,
$-\log(1+x)\geq -\log(1+\bar x)-\frac{x-\bar x}{1+\bar x}$, and
$x^2\geq  
2 \bar xx-\bar x^2$,
for any $\bar x,\bar y$ in the function domains.
Thus, given a feasible local point $(\ab^r[n],\qb^r[n])$, the rate function \eqref{ratefunction} admits a global lower bound, shown as follows
\begin{subequations} \label{eqn: lower bound}
\begin{align}
&\!\!R_k(\ab[n],\qb[n])\notag\\
&= \log\bigg(1+  \sum_{j=1}^K     \frac{\gamma a_j^2[n]}{d_{jk}[n]} \bigg)
-\log\bigg(1+  \sum_{j\neq k}^K \frac{\gamma a_j^2[n]}{d_{jk}[n]} \bigg)\\
&\geq\log\!\bigg(\!1\!+\!\gamma\sum_{j=1}^K \bigg[\frac{2 a_j^r[n]}{d_{jk}^r[n]}{a_j[n]} - \frac{p_j^r[n]}{(d_{jk}^r[n])^2}d_{jk}[n] \bigg]  \bigg)\notag\\
&~~~~~\!\!-\!\log(1\!+\!I_k^r[n]) \!- \!\frac{I_k[n]\!-\!I_k^r[n]}{1 + I_k^r[n]}
 \label{eqn: lower bound e1} \\
&\geq\log\bigg(1\!+\!\gamma\sum_{j=1}^K \bigg[\frac{2 a_j^r[n]}{d_{jk}^r[n]}{a_j[n]} \!-\! \frac{p_j^r[n]}{(d_{jk}^r[n])^2}\|\qb_j[n]-\sb_k\|^2 \bigg]  \bigg)
  \notag\\
&~~~~ -\log(1\!+\! I_k^r[n])+\frac{I_k^r[n]}{1 + I_k^r[n]}\notag\\
&~~~~-\!\frac{\gamma}{1 + I_k^r[n]}\sum_{j\neq k}^K \frac{ a_j^2[n]}{d_{jk}^r[n] \!+\! 2(\qb_j^r[n]\!-\!\sb_k)^T(\qb_j[n]\!-\!\qb_j^r[n])}   \label{eqn: lower bound e2} \\
&\triangleq \tilde R_k(\ab[n],\qb[n]; \ab^r[n], \qb^r[n]), 
\end{align}
\end{subequations}
where { \begin{align}\smaller[0.5] I_k^r[n]=\sum_{j\neq k}^K \frac{\gamma (a_j^r[n])^2}{d_{jk}^r[n]},\notag \end{align}} 
\!\!inequality in \eqref{eqn: lower bound e1} is obtained by the linear lower bounds of the convex function $\frac{\alpha_j^2[n]}{d_{jk}[n]}$ and the logarithm function $-\log(1+I_k[n])$, while \eqref{eqn: lower bound e2} is obtained by the linear lower bound of the convex function $\| \qb_j[n]-\sb_k\|^2$.
The lower bound $\tilde R_k(\ab[n],\qb[n]; \ab^r[n], \qb^r[n])$ is a concave function in $(\ab[n],\qb[n])$, and it is locally tight, i.e.,
\begin{align}\label{eqn: local tight}
R_k(\ab^r[n], \qb^r[n]) = \tilde R_k(\ab^r[n], \qb^r[n]; \ab^r[n], \qb^r[n]).
\end{align}
Similarly, to handle the non-convex constraint \eqref{eqn:CA constraint}, we apply the convex lower bound of $\| \qb_k[n]-\qb_j[n]\|^2$, which gives
\begin{align}\label{eq: CA approx}
2(\qb^r_k[n]-  \qb^r_j[n])^T{(\qb_k[n]-\qb_j[n])} \geq  (d_{jk}^r[n])^2 + d_{\min}^2.
\end{align}

Therefore, given $\{\qb^r[n], \ab^r[n]\}_{n=0}^{M+1}$ at the $r$th iteration, the proposed SCA method solves the following convex approximation problem for \eqref{Problem::SEM3}
\begin{subequations}\label{Problem::SEM:SCA}
	\begin{align}
	&\!\!\!\!\!\!\!\!\!\!\!\{\qb^{r+1}[n],\ab^{r+1}[n]\}_{n=1}^{M}=\notag\\
	& \argmax_{\{\ab[n],\, \qb[n]\}_{n=1}^{M}} ~\sum_{n=1}^M\sum_{k=1}^{K} \tilde R_k(\ab[n],\qb[n];  \ab^r[n],  \qb^r[n])\notag\\
	&~~~\st
	~	\eqref{altitute constraint}, \eqref{FullSpeedLimit}, \eqref{eq: CA approx} {~\rm for~}n\in \Mset, j,k\in \Kset,j>k, \\
	&~~~~~~~~ 0 \leq  a_k[n]\leq \sqrt{P_{\max}},~\forall k\in \Kset, n\in \Mset. \label{eqn::SEM SCA C1}
	\end{align}
\end{subequations}	
In Algorithm \ref{Alg1}, we summarize the proposed SCA-based TPC algorithm for solving \eqref{Problem::SEM3}.
\begin{algorithm}
	\caption{SCA-based TPC algorithm for solving \eqref{Problem::SEM3}.} \label{Alg1}\smaller[1]
	\begin{algorithmic}[1]
		\STATE {Set} $r=0$, and iteration tolerance $\epsilon>0$.
		\STATE {Initialize} the trajectory $\qb_k^{0}[n]$ and power vector $p_k^{0}[n]$ for  $k\in\Kset$, $n\in \Mset$.
		\STATE {Obtain} $R^0 = \sum_{k,n}R_k(\pb^0[n], \qb^{0}[n])$ by \eqref{Rate}.
		\REPEAT
		\STATE Update $\{a_k^r[n],d_{jk}^r[n],I_k^r[n]\}$ with $\{\qb_k^r[n],p_k^r[n]\}$ by \eqref{Eq11} and \eqref{Eq12}.
		\STATE Update $\{a_k^{r+1}[n],\qb_k^{r+1}[n]\}$ by solving the convex optimization problem \eqref{Problem::SEM:SCA}.
		\STATE Update $R^{r+1} = \sum_{k=1}^K \sum_{n=1}^MR_k(\ab^{r+1}[n], \qb^{r+1}[n])$ by \eqref{Rate}.
		\STATE Set $r:=r+1$.		
		\UNTIL $\frac{R^{r} - R^{r-1}}{R^{r-1}}\leq \epsilon$.
		\STATE \textbf{Output} $\{p_k^{r}[n],\qb_k^{r}[n]\}$
	\end{algorithmic}
\end{algorithm}

\begin{rem}\rm {\bf(Convergence)} It can be verified that the  optimal objective value of \eqref{Problem::SEM:SCA} is non-decreasing with the iteration number $r$ and it is upper-bounded. Therefore, the sequence of objective values of \eqref{Problem::SEM:SCA} converge as $r$ goes to infinity.
\end{rem}

\begin{rem}\rm  {\bf(Comparison with AO)} It is possible to solve problem \eqref{Problem::SEM} by the AO technique \cite{T:Relay:1S1R1Dzeng,Wu2018TWC}, which optimizes the aggregate sum rate with respect to $\{\qb[n]\}$ and $\{\pb[n]\}$ in an alternating manner. When the trajectories $\{\qb[n]\}$ are fixed, the update of the powers  $\{\pb[n]\}$ can be achieved, e.g., by the WMMSE approach \cite{WMMSE}. When the powers $\{\pb[n]\}$  are fixed, the update of the trajectories $\{\qb[n]\}$ may be achieved by the SCA technique similar to Algorithm \ref{Alg1}. It will be shown in Section \ref{sec:simulation} that the AO method can yield comparable aggregate sum rate as the Algorithm \ref{Alg1}, but it is computationally more expensive.
\end{rem}

%
%

\section{Time-Efficient TPC Optimization}\label{sec: distr opt}
In contrast to the original problem \eqref{Problem::SEM} which involves $\mathcal{O}(KN)$  variables and constraints, the order of variables and constraints of \eqref{Problem::SEM3} is $\mathcal{O}(KM)$,  which is usually  much smaller. 
However, in practical scenarios, the number of UAVs $K$ and/or the number of time slots $M$ can be large (e.g., when the UAVs' initial locations are far away from their serving GTs' locations). Thus, problem \eqref{Problem::SEM3} can still be computationally expensive  and it may require a long computation time. 

In this section, we aim to overcome this computational  issue from two different perspectives. Firstly, in Section \ref{subsec: distributed alg}, we propose a parallel algorithm for solving the TPC problem \eqref{Problem::SEM3}.
Such a parallel algorithm can be implemented over multi-core computers/computer clusters, which thus has significantly reduced computation time than Algorithm \ref{Alg1}.
Secondly, in Section \ref{subsec: waypoint model}, we propose a segment-by-segment method which decomposes the TPC problem for the time interval $[0,M]$ into successive TPC problems each of which involves only a smaller flight interval.

\subsection{Parallel TPC Optimization}\label{subsec: distributed alg}

The proposed parallel TPC algorithm is based on the BSUMM method \cite{BSUM}. The BSUMM method is a generalization of the classical alternating direction method of multipliers (ADMM) \cite{ADMM}, in the sense that each block updating maximizes a locally tight surrogate function rather than the objective function of the original problem as in ADMM.
As will be shown shortly, the locally tight surrogate function for problem \eqref{Problem::SEM3} is carefully designed so that the variable updates can be performed in parallel, each for one UAV.

To present the proposed parallel algorithm, let us define
\begin{align}\label{eqn: zjk}
       \zb_{kj}[n] =\qb_k[n]-\qb_j[n],~\forall k,j\in \Kset,~ j>k,~n\in\Mset,
\end{align}
and let $\zb[n]\in\Rset^{3K(K-1)}$ be a vector concatenating all the $\zb_{kj}[n]$, $k,j\in \Kset,~ j>k$. 
Let $\bar\Ab \in \mathbb{R}^{\tfrac{1}{2}K(K-1)\times K}$ be a matrix obtained by vertically stacking $[ {\bf 0}_{K-k,k-1},{\bf 1}_{K-k,1},-\Ib_{K-k} ]$ for $k=1,\ldots,K-1$, and let
$\Ab=\bar\Ab \otimes \Ib_3$. 
Then \eqref{eqn: zjk} can be compactly expressed as
\begin{align}
   \Ab\qb[n]=\zb[n],~n\in \Mset.
\end{align}
Besides, 
let us define the feasible set of variables as
\begin{align}
\!\!\!\Zset  &\!\triangleq\! \big\{ \zb\in\Rset^3~|~\|\zb \|\geq\dmin\big\},\\
\!\!\!\Xset_k&\!\triangleq\!  \{(a_k[n], \qb_k[n]),n\!\in\! \Mset|  \eqref{altitute constraint}, \eqref{FullSpeedLimit}, \eqref{eqn::SEM SCA C1}, \qb_k[M]\!=\!\qb_k^*\},
\end{align}
for all $k\in \Kset$. Then, we can write problem \eqref{Problem::SEM3} as
\begin{subequations}\label{Problem::SEM4}
	\begin{align}	
	\!\!\!\max_{\{\ab[n],\qb[n], \zb[n]\}}&\sum_{n=1}^M \sum_{k=1}^{K}  R_k(\ab[n],\qb[n])  \\
	\!\!\!\st~~~~~&
	\{(a_k[n], \qb_k[n])\}_{n=1}^M\in \Xset_k,~\forall k\in \Kset, \\
	\!\!\!&\zb_{kj}[n] \in \Zset,~\forall k,j\in \Kset,~ j>k,~n\in\Mset, \\
	\!\!\!&  \Ab\qb[n]=\zb[n],~\forall n\in \Mset. \label{eqn::SEM eq form C1}
	\end{align}
\end{subequations}

Let us first apply the standard ADMM \cite{ADMM} to \eqref{Problem::SEM4}. Let $\lambdab_{kj}[n]\in\Rset^3$ be the Lagrange dual variable associated with each linear constraint $\zb_{kj}[n]=\qb_k[n]-\qb_j[n]$ in  \eqref{eqn::SEM eq form C1} and $\lambdab[n]$ be a vector concatenating all $\lambdab_{ij}[n]$. Then the partial augmented Lagrangian of \eqref{Problem::SEM4} is given by
\begin{align}
\mathcal{L} =\sum_{n=1}^M\bigg(\sum_{k=1}^{K}  R_k(\ab[n],\qb[n]) - \lambdab^T[n](\Ab\qb[n] - \zb[n]) - \frac{1}{2}\|\Ab\qb[n]-\zb[n]\|_\Bb^2\bigg),
\end{align}
where $\Bb=\diag(\bb\otimes{\bf 1}_3)$ is a diagonal matrix with $\bb = [b_{kj}]_{k,j\neq k}\in\Rset_{++}^{K(K-1)}$ being some penalty parameters.
By the standard ADMM, we have the following iterative updates for problem \eqref{Problem::SEM4},
\hspace{-0.5cm}
\begin{subequations}
\hspace{-0.5cm}\begin{align} 
\hspace{-0.5cm}\!\!\!\!\!\!\!
&\{\ab^{r+1}[n], \qb^{r+1}[n]\}_{n=1}^M \notag \\
& =\argmax_{\substack {\ab[n],\qb[n],\\ n\in\Mset}} 
\bigg\{\sum_{n=1}^M \sum_{k=1}^{K}  R_k(\ab[n],\qb[n]) 
-\frac{1}{2}\sum_{n=1}^M\|\Ab\qb[n]-\zb^r[n] + \Bb^{-1}\lambdab^r[n]\|^2_{\Bb}\bigg\},   \notag  \\ 
&\qquad\qquad\quad\st~~\{(a_k[n], \qb_k[n])\}_{n=1}^M\in \Xset_k,~k\in \Kset,  \label{eqn: admm update x} 
\\
%
&\zb_{kj}^{r+1} [n]  \!=\! \argmin_{\zb_{kj}[n]\in\Zset}\! \bigg\{\frac{b_{kj}}{2}\|\qb_k^{r\!+\!1}[n] -\! \qb_j^{r\!+\!1}[n]  \! \! -\! \zb_{kj}[n] \!+\! b_{kj}^{-1} \lambdab_{kj}^r[n] \|_2^2\bigg\},  \label{eqn: admm update z}  
\\
&\lambdab_{kj}^{r+1} [n]\!=\! \lambdab_{kj}^{r} [n] \!+\! b_{kj} (\qb_k^{r+1}[n] \!-\! \qb_j^{r+1}[n]  \!-\! \zb_{kj}^{r+1}[n]), \label{eqn:admm update lambda} 
\end{align} 
\end{subequations}
where $n\!\in\!\Mset$, $k,j\in \Kset,j>k$ in \eqref{eqn: admm update z} and \eqref{eqn:admm update lambda}, and the superscript $r$ denotes the iteration index.

Note that although the feasible set $\Zset$ is non-convex, the update of $\zb[n]$ in \eqref{eqn: admm update z} (i.e., projecting $\qb_k^{r+1}[n] - \qb_j^{r+1}[n] + b_{kj}^{-1} \lambdab_{kj}^r[n]$ onto $\Zset$) has a closed-form solution as
\begin{align}
\zb_{kj}^{r+1} [n] =  \frac{ \qb_k^{r+1}[n] - \qb_j^{r+1}[n] +  b_{kj}^{-1} \lambdab_{kj}^r[n]  }{\min\left\{  \frac{\| \qb_k^{r+1}[n] - \qb_j^{r+1}[n]  +  b_{kj}^{-1} \lambdab_{kj}^r[n]\|}{d_{\min}} ,1 \right\} }.
\label{eqn:admm update z solution}
\end{align} Therefore, both updates in \eqref{eqn: admm update z} and \eqref{eqn:admm update lambda} are simple and separable with respect to the UAVs. 

However, \eqref{eqn: admm update x} is non-convex and difficult to handle. Moreover, the objective function of \eqref{eqn: admm update x} couples all the optimization variables.
To resolve this issue, we adopt the BSUMM strategy to derive a concave, separable lower bound for the objective function of \eqref{eqn: admm update x} that is amenable to efficient distributed and parallel updates.
Specifically, the derived concave lower bound is given by
\begin{multline}\label{eqn:surrogate func}
\sum_{k=1}^{K} \sum_{n=1}^M\hat R_k(a_k[n],\qb_k[n]; \ab^r[n], \qb^r[n])\\[-8pt]
 -\frac{1}{2}\sum_{n=1}^M\|\Ab\qb[n]-\zb^r[n] + \Bb^{-1}\lambdab^r[n]\|^2_{\Bb}\\[-8pt]
 - \frac{1}{2}\sum_{n=1}^M \|\qb[n] - \qb^r[n]\|^2_{\Cb-\Ab^T\Bb\Ab}.
\end{multline}  

The first term $ \sum_{k=1}^K \sum_{n=1}^M \hat R_k(a_k[n],\qb_k[n]; \ab^r[n], \qb^r[n])$ will be shown to be a concave lower bound of the rate function $\sum_{n=1}^M \sum_{k=1}^{K}  R_k(\pb[n],\qb[n])$;
the third term in \eqref{eqn:surrogate func} is a proximal penalty with a diagonal matrix $\Cb$ satisfying $\Cb-\Ab^T\Bb\Ab \succeq \zerob$.

To obtain the first term in \eqref{eqn:surrogate func}, denote
\begin{align}\label{eqn: coefficient}
\mu_{jk}^r[n] = \frac{ p_j^r[n]/d_{jk}^r[n]+\varepsilon} 
{\sum_{i=1}^K \frac{p_i^r[n]}{d_{ik}^r[n]} +K\varepsilon},~\forall k,j\in \Kset,
\end{align}
where $\varepsilon$ is a small positive number such that $\mu_{jk}^r[n]>0$. 
Note that $\sum_{j=1}^K \mu_{jk}^r[n] =1$.
By  \eqref{eqn: lower bound} and \eqref{eqn: coefficient}, we have a lower bound for the term $\sum_{k=1}^{K}  R_k(\ab[n],\qb[n])$ in \eqref{eqn: admm update x}, which
is given by 
\begin{align}
&\sum_{k=1}^KR_k(\ab[n],\qb[n])  \geq \sum_{k=1}^K \tilde{R}_k(\ab[n],\qb[n]; \ab^r[n],\qb^r[n])  \notag \\[-1pt]
&=\sum_{k=1}^K\log\!\bigg(1\!+\!\gamma \sum_{j=1}^K \!\frac{\mu_{jk}^r[n]}{\mu_{jk}^r[n]} \bigg[ \frac{2a_j^r[n]a_j[n]}{d_{jk}^r[n]} - \frac{p_j^r[n] }{(d_{jk}^r[n])^2} \|\qb_j[n]-\sb_k\|^2\bigg]\bigg) -\sum_{k=1}^K \log(1+ I_k^r[n])\notag \\[-1pt]
&~~~~~~~~~~~~~~~~~~~~~~~~~ + \sum_{k=1}^K\frac{I_k^r[n]}{1 + I_k^r[n]} - \sum_{k=1}^K\frac{\gamma}{1+I_k^r[n]}
\sum_{j\neq k}^K \frac{a_j^2[n]}{d_{jk}^r[n] + 2(\qb_j^r[n]-\sb_k)^T( \qb_j[n]-\qb_j^r[n]) } \notag\\[-1pt]
&\geq\sum_{k=1}^K\sum_{j=1}^K \mu_{jk}^r[n] \log\bigg(1+ \frac{\gamma}{\mu_{jk}^r[n]}\bigg[\frac{2a_j^r[n]a_j[n]}{d_{jk}^r[n]} - \frac{p_j^r[n] }{(d_{jk}^r[n])^2} \|\qb_j[n]-\sb_k\|^2\bigg]  \bigg) -  \sum_{k=1}^K\log(1+ I_k^r[n]) \notag\\
&~~~~~~~~~~~~~~~~~~~~~~~~~ +\sum_{k=1}^K\frac{I_k^r[n]}{1 + I_k^r[n]} 
- \sum_{k=1}^K\frac{\gamma}{1+I_k^r[n]}\sum_{j\neq k}^K\frac{ \alpha_j^2[n]}{d_{jk}^r[n] + 2(\qb_j^r[n]-\sb_k)^T( \qb_j[n]-\qb_j^r[n]) }\notag
\end{align}
\begin{align}
&= \sum_{k=1}^K\Bigg\{\sum_{j=1}^K \mu_{kj}^r[n] \log\bigg(  1+\frac{ \gamma}{\mu_{kj}^r[n]} \bigg[ \frac{2a_k^r[n]a_k[n]}{d_{kj}^r[n]} - \frac{p_k^r[n] }{(d_{kj}^r[n])^2} \|\qb_k[n]-\sb_j\|^2\bigg]  \bigg) -  \log(1+ I_k^r[n])  \notag \\[-1pt]
&~~~~~~~~~~~~~~~~~~~~~~~~~~~~~~~~~ +\frac{I_k^r[n]}{1 + I_k^r[n]} 
- \sum_{j\neq k}^K\frac{\gamma}{1+ I_j^r[n]} \frac{a_k^2[n]}{d_{kj}^r[n] + 2(\qb_k^r[n]-\sb_j)^T( \qb_k[n]-\qb_k^r[n]) }\Bigg\} \notag \\[-1pt]
& \triangleq \sum_{k=1}^K\hat R_k(a_k[n],\qb_k[n]; \ab^r[n], \qb^r[n]). \label{eqn: lower bound 2}
\end{align} 

In \eqref{eqn: lower bound 2}, the second inequality is due to the concave logarithm function.
Notice that each $\hat R_k(a_k[n],\qb_k[n]; \\ \ab^r[n], \qb^r[n])$ involves only variables $ \{a_k[n], \qb_k[n]\}_{n=1}^M$ associated with UAV $k$. Therefore, the lower bound function in \eqref{eqn: lower bound 2} is decomposable across $K$ UAVs. Besides, when $\varepsilon\rightarrow 0$, the function is locally tight, i.e., $\sum_{k=1}^K R_k(\ab^r[n],\qb^r[n])=\sum_{k=1}^K \hat R_k(a_k^r[n],\qb_k^r[n]; \ab^r[n], \qb^r[n]).$

It can be shown that the second and third terms in \eqref{eqn:surrogate func} are also decomposable across $K$ UAVs.
In particular, let $\Cb=\diag(\cb\otimes{\bf 1}_3)$ with $\cb=[c_1,\ldots,c_K]^T\in\Rset_{++}^K$. Then, we have
\begin{align} 
&\!\!\sum_{n=1}^M\!\|\Ab\qb[n]\!-\!\zb^r\![n] \!+\! \Bb^{-1}\!\lambdab^r[n] \|^2_{\Bb} \!+\! \sum_{n=1}^M\! \|\qb[n] \!-\! \qb^r\![n]\|^2_{\Cb\!-\!\Ab^T\!\Bb\Ab} \notag \\[-1pt]
&\equiv \sum_{n=1}^M \bigg(\| \qb[n] - \qb^r[n] + \Cb^{-1}\Ab^T\Bb(\Ab\qb^r[n]-\zb^r[n]+\Bb^{-1}\lambdab^r[n])\|^2_{\Cb}\bigg) 
\notag 
\\
%
&=\sum_{k=1}^K\bigg ({c_k}\!\sum_{n=1}^M \big\| \qb_k[n] - \qb_k^r[n] +c_k^{-1}\Ab_k^T\Bb(\Ab\qb^r[n]-\zb^r[n]+\Bb^{-1}\lambdab^r[n])\big\|_2^2\bigg)\notag\\[-3pt]
&=\sum_{k=1}^K {c_k} \sum_{n=1}^M \| \qb_k[n] - \hat\qb_k^r[n]\|_2^2,  \label{eqn: linearization}
\end{align}
where `$\equiv$' means equivalence up to a constant, $\Ab_k\in\Rset^{3K(K-1)\times 3}$ is the $k$th column block of $\Ab$, i.e., $\Ab=[\Ab_1,\ldots,\Ab_K]$, and 
$\hat\qb_k^r[n] \triangleq \qb_k^r[n] - c_k^{-1}\Ab_k^T\Bb(\Ab\qb^r[n] - \zb^r[n]+\Bb^{-1}\lambdab^r[n]).$

By replacing the objective function of  \eqref{eqn: admm update x} by \eqref{eqn:surrogate func} and \eqref{eqn: linearization}, we have the following $K$ parallel convex subproblems 
\begin{align}\label{eqn: admm update x separate}
&\!\!\!\!\!\!\{a_k^{r+1}[n], \qb_k^{r+1}[n]\}_{n=1}^M= \notag \\[-2pt]
&\argmax_{\substack { a_k[n],\qb_k[n],\\ n\in\Mset }} 
\bigg\{\sum_{n=1}^M \Big(\hat R_k(a_k[n],\qb_k[n]; \ab^r[n], \qb^r[n]) - \frac{c_k}{2}\| \qb_k[n] - \hat\qb_k^r[n]\|^2\Big)\bigg\} \\
&~~~~~\st~~~\{(a_k[n], \qb_k[n])\}_{n=1}^M\in \Xset_k, \notag
\end{align}
for all $k\in \Kset.$
The proposed parallel TPC algorithm is summarized in Algorithm \ref{table: ADMM for SEM}.
We should emphasize that the proposed Algorithm \ref{Alg2} enables parallel computation over multi-core CPUs and therefore can greatly reduce the computation time, as will be demonstrated in Section \ref{sec:simulation}. \vspace{-.1\baselineskip}
\begin{algorithm}[!b]\centering
	\caption{Proposed parallel TPC algorithm for \eqref{Problem::SEM3}. }\label{Alg2}\smaller[1]
	\begin{algorithmic}[1]\label{table: ADMM for SEM}
		\STATE {\bf Set} $r=0$, $\bb\succeq{\bf 0}$, $\cb\succeq{\bf 0}$ such that $\Cb-\Ab^T\Bb\Ab\succeq{\bf 0}$.
		\STATE {\bf Initialize} $\lambdab_{kj}^0[n]$, $a_k^0[n]$ and $\qb_k^0[n]~\forall k,j\in \Kset$ and $n\in \Mset$.\!\!
		\REPEAT
		\FORALL{$k=1,\ldots,K$ (in parallel)}
		\STATE Update $\hat\qb_k^r[n] \!=\!\qb_k^r[n] \!-\! c_k^{-1}\Ab_k^T\Bb(\Ab\qb^r[n] \!-\! \zb^r[n]+\Bb^{-1}\lambdab^r[n])$\!\!
		\STATE Update $\{a^{r+1}_k[n], \qb^{r+1}_k[n]\}_{n=1}^N$ by solving \eqref{eqn: admm update x separate}.
		\STATE Update $\zb_{kj}^{r+1}[n]$ by \eqref{eqn:admm update z solution} for all $j\!\in\!\Kset,~j\!>\!k,$, $n\!\in\!\Mset$.\!\!
		\STATE Update $\lambdab_{kj}^{r+1} [n]$ by \eqref{eqn:admm update lambda} for all $j\!\in\!\Kset,~j>k$, $n\!\in\!\Mset$.		
		\ENDFOR 		
		\STATE {\bf Set} $r =r+1$.
		\UNTIL predefined stopping condition is satisfied.
	\end{algorithmic}
\end{algorithm}\vspace{-5pt}

\subsection{Segment-by-Segment TPC Optimization}\label{subsec: waypoint model}
The proposed parallel TPC algorithm decomposes the problem into $K$ parallel subproblems in each iteration. However, when the number of discrete time slots $M$ is large, problem \eqref{Problem::SEM3} is still time-consuming to solve. To overcome this issue, a simple idea is to divide the time interval $[0,M]$ into several smaller time segments and apply the (parallel) TPC algorithms to each of them.

Specifically, denote $N_{\rm seg}$ as the number of time slots per segment. The time slot set of the $\ell$th segment is given by
\begin{align}
 \Nset_\ell = \{(\ell-1)N_{\rm seg}+1, (\ell-1)N_{\rm seg}+2, \ldots, \ell N_{\rm seg} \}.
\end{align}
%
Given $\{\qb_k[n],p_k[n],\forall k\}_{n\in\Nset_{\ell-1}}$ and similar to \eqref{Problem::SEM3},
we determine the trajectories and powers for the $\ell$th segment by solving the following problem:\vspace{-5pt}
\begin{subequations}\label{Problem::SbyS}
	\begin{align}	
	&\!\!\!\!\!\{\pb^\star[n],\qb^\star[n]\}_{n\in\Nset_{\ell}} \!=\! \argmax_{\substack {p_k[n],\qb_k[n]\\ k\in\Kset,n\in \Nset_{\ell}} }\sum_{n\in \Nset_{\ell}} \sum_{k=1}^{K}  R_k(\pb[n],\qb[n])\\
	&\hspace{15mm} \st~\eqref{altitute constraint}, \eqref{FullSpeedLimit}, \eqref{eqn:CA constraint},\eqref{eqn::SEM C1}~\forall n\in  \Nset_\ell, k\in \Kset, \\
	&\hspace{15mm}~~~~~\qb[(\ell-1)N_{\rm seg}]=\qb^\star[(\ell-1)N_{\rm seg}], 
	\end{align}
\end{subequations}
where $\qb^\star[0]=\qb[0]$. 
Note that the segment-by-segment TPC problem \eqref{Problem::SbyS} is solved sequentially for each $\ell$ until the sum rate $\sum_{k=1}^{K}  R_k(\pb^\star[ \ell N_{\rm seg}],\qb^\star[ \ell N_{\rm seg}])$ achieves at least the same value as that obtained by \eqref{Problem::DP}. 
We summarize the segment-by-segment method in Algorithm \ref{Alg3}.
\begin{algorithm}[!t]\centering
	\caption{Proposed segment-by-segment TPC algorithm for \eqref{Problem::SEM3}. }\label{Alg3}\smaller[1]
	\begin{algorithmic}[1]
		\STATE {\bf Given} $(\qb[0],\pb[0], R_s[0])$ and a positive integer $N_{\rm seg}$. Set $\ell=0$.
		\STATE {\bf Find} the maximum sum rate $R_s^\star$ when all UAVs are in the hovering state by solving \eqref{Problem::DP}.
		
		\WHILE{$R_s[(\ell-1)N_{\rm seg}]<R_s^\star$}
		\STATE Obtain $\{\qb^\star[n],\pb^\star[n]\}_{n\in \Nset_\ell}$ by solving \eqref{Problem::SbyS} with Algorithm \ref{Alg2}.
		\STATE $\ell:=\ell + 1$. 
		\ENDWHILE
		
		\STATE \textbf{Output} the TPC solution $\{\qb^\star[n],\pb^\star[n]\}$.
	\end{algorithmic}
\end{algorithm}\vspace{-5pt}

\section{TPC Optimization over Orthogonal Channels}\label{sec: FDMA/TDMA}

The $K$-user UAV system may operate under the FDMA or the TDMA schemes, which is particularly interesting when 
the interference channels are strong \cite{tse2005fundamentals}. This section extends the TPC optimization methods to the FDMA/TDMA schemes.

Denote $\alpha_k[n]\in [0,1]$ as the fraction of time/bandwidth resource allocated to the UAV-GT pair $k$ at time slot $n$, which satisfies $\sum_{k=1}^K\alpha_k[n] = 1$ for all $n\in \Nset$.
The achievable rates in bps/Hz of the $k$th UAV-GT pair at time slot $n$ under FDMA and TDMA schemes are respectively given by 
\begin{align}
R_k^{\rm FDMA}(\alpha_k[n],\qb_k[n])&\!\triangleq\! \alpha_k[n] \log_2\!\bigg(\! 1\!+\! \frac{\gamma \Pmax}{\alpha_k[n]\|\qb_k[n]\!-\!\sb_k\|^2}\!\!\bigg),\\
R_k^{\rm TDMA}(\alpha_k[n],\qb_k[n])&\!\triangleq\! \alpha_k[n] \log_2\!\bigg(\! 1\!+\! \frac{\gamma \Pmax}{           \|\qb_k[n]\!-\!\sb_k\|^2}\!\bigg).\!\!
\end{align} 
Here we only consider the short-term peak power constraints for the UAV transmission.
It is easy to see that FDMA does not perform worse than the TDMA since $R_k^{\rm FDMA}(\alpha_k[n],\qb_k[n]) \geq R_k^{\rm TDMA}(\alpha_k[n],\qb_k[n])$ for all $k\in \Kset$ and $n\in \Nset.$
The trajectory and resource allocation problem for FDMA/TDMA is given by 
\begin{subequations}\label{Problem::SEM:OMA}
	\begin{align}	
	\max_{\{\alpha_k[n],\qb_k[n]\}}&~\sum_{k=1}^K\sum_{n=1}^N  R_k^{\rm XDMA}(\alpha_k[n],\qb_k[n])  \label{E8:a}\\
	\st~~	
	& \eqref{altitute constraint}, \eqref{FullSpeedLimit}, \eqref{eqn:CA constraint} {~\rm for~}n\in  \Nset, k\in \Kset, \\
	& 0\leq \alpha_k[n]\leq 1,~\sum_{k=1}^K \alpha_k[n] = 1,~\forall k,n,\label{E8:f}
	\end{align}
\end{subequations}	
where XDMA refers to TDMA or FDMA. 


Problem \eqref{Problem::SEM:OMA} can be handled in a similar fashion as Algorithm \ref{Alg1} based on the SCA technique.
Specifically, a concave lower bound for $R_k^{\rm FDMA}(\alpha_k[n],\qb_k[n])$ can be derived as
\begin{align}\label{FDMA:LB}
	&R_k^{\rm FDMA}(\alpha_k[n],\qb_k[n])\notag\\
	&  =\alpha_k[n] \log_2\!\bigg(\!  1\!+\! \frac{\gamma \Pmax}{\alpha_k[n] \|\qb_k[n]-\sb_k\|^2}\!\bigg) \notag \\
	&\ge \alpha_k[n] \log_2\!\bigg(\! 1\!+\! \frac{\gamma \Pmax}{\alpha_k[n]}\!
						   \bigg[\!\frac{3}{\|\qb_k^r[n]\!-\!\sb_k\|^2} \!-\! \frac{2(\|\qb_k[n]\!-\!\sb_k\|)}{\|\qb_k^r[n]\!-\!\sb_k\|^3}\!\bigg]\!\bigg) \vspace{-10pt}
\end{align}
which is concave in both $\alpha_k[n]$ and $\qb_k[n]$.

For the TDMA scheme, firstly, note that problem \eqref{Problem::SEM:OMA} admits an optimal time allocation 
\begin{subnumcases}{\!\!\!\alpha_k^*[n] = \label{TDMA:OptimalTA}}
\!\!1, & \!\!\!\!\!\!if $k=\argmin_{j\in\Kset}\|\qb_j[n]-\sb_j\|$,\\
\!\!0, & \!\!\!\!\!\!otherwise,
\end{subnumcases}
for all $n\in\Nset$. Under \eqref{TDMA:OptimalTA}, we have
\begin{subequations}
	\begin{align}
	&\sum_{k=1}^K R_k^{\rm TDMA}(\alpha_k[n],\qb_k[n])\\[-2pt]
   &=\sum_{k=1}^K \alpha_k[n] \log_2\bigg( 1+ \frac{\gamma \Pmax}{\|\qb_k[n]-\sb_k\|^2}\bigg)\\[-2pt]
   &=\log_2\!\bigg(\! 1\!+\! \sum_{k=1}^K \!\frac{\gamma \Pmax \alpha_k[n]}{\|\qb_k[n]\!-\!\sb_k\|^2}\bigg) 
   \!\!=\! \log_2\!\bigg(\! 1\!+\! \sum_{k=1}^K \!\frac{\gamma \Pmax \beta^2_k[n]}{\|\qb_k[n]\!-\!\sb_k\|^2}\bigg)\notag
   \\
&\geq\log_2\bigg( \!1\!+\! \sum_{k=1}^K \bigg[\frac{ 2\gamma \Pmax\beta_k^r[n]\beta_k[n]}{\|\qb_k^r[n]-\sb_k\|^2}  - \frac{\gamma \Pmax (\beta_k^r)^2[n]}{\|\qb_k^r[n]-\sb_k\|^4} (\|\qb_k[n]-\sb_k\|^2)\bigg]\bigg),	  \label{E11c}
	\end{align}
\end{subequations}
where we have defined $\beta_k[n]\triangleq \sqrt{\alpha_k[n]}$  and applied the linear lower bound as used in  \eqref{eqn: lower bound}.
The function in \eqref{E11c} is a locally tight concave lower bound of $\sum_{k=1}^K R_k^{\rm TDMA}(\alpha_k[n],\qb_k[n])$.

\section{Simulation Results}\label{sec:simulation}
\subsection{Simulation Settings}
Unless otherwise stated, all GTs are randomly located in an area of 1 square kilometer centered at $(0,0,0)$, and the altitude of each UAV is initialized with $\Hmin\!=\!100$ m. All UAVs have the same flying time $T\!=\!10$ minutes, maximum power limit $\Pmax\!=\!30$ dBm, and speed limits $V_L\!=\!$ 20 m/s, $V_A\!=\!$ 5 m/s, and  $V_D\!=\!$ 3 m/s. We set $B=10$ MHz, $\beta_0\!=\!-50$ dB, and $N_0\!=\!-160$ dBm/Hz.
By \eqref{maxTs}, we set $T_s\!=\! \frac{\dmin}{\sqrt{4V_L^2 + (V_D + V_A)^2}}$, which is about 0.4903 seconds for the case of $\dmin\!=\!20$ m. For the proposed parallel TPC Algorithm \ref{Alg2}, we set $b_{kj}\!=\!b\!=\!0.001$ and $c_k \!=\! c \!=\! 1.1b\lambda_{\max}(\Ab^T\!\Ab)$ for all $k$ and $j$, where $\lambda_{\max}(\cdot)$ denotes the maximum eigenvalue. 
A simple pre-conditioning is performed where the considered optimization problems are scaled such that $\Hmin$ is normalized to be 1.

All convex problems (e.g. \eqref{Problem::SEM:SCA}, \eqref{eqn: admm update x separate}) are solved by a customized logarithmic barrier based interior point method (IPM) (specifically, see \cite[Algorithm 11.1]{boyd2004convex}) with barrier parameters $s=1$ and $\mu=30$. 
The centering problems are solved using the Newton's method  \cite[Algorithm 10.1]{boyd2004convex} and the associated step size is chosen by the backtracking line search \cite[Algorithm 9.2]{boyd2004convex} with parameters $\alpha = 0.01$ and $\beta=0.5$.
The inner iteration of the Newton update is terminated if the iteration number is larger than 30 or a relative precision $10^{-8}$ is satisfied, while the outer IPM iteration stops whenever $M/s \leq 10^{-8}$. 

\subsection{Impact of Flight Altitude}

\begin{figure}\centering
	\subfigure[Initial trajectories on x-y plane  with fixed height.]
			  {\label{fig:InitTraj:FixedH}\includegraphics[width=0.4\linewidth]{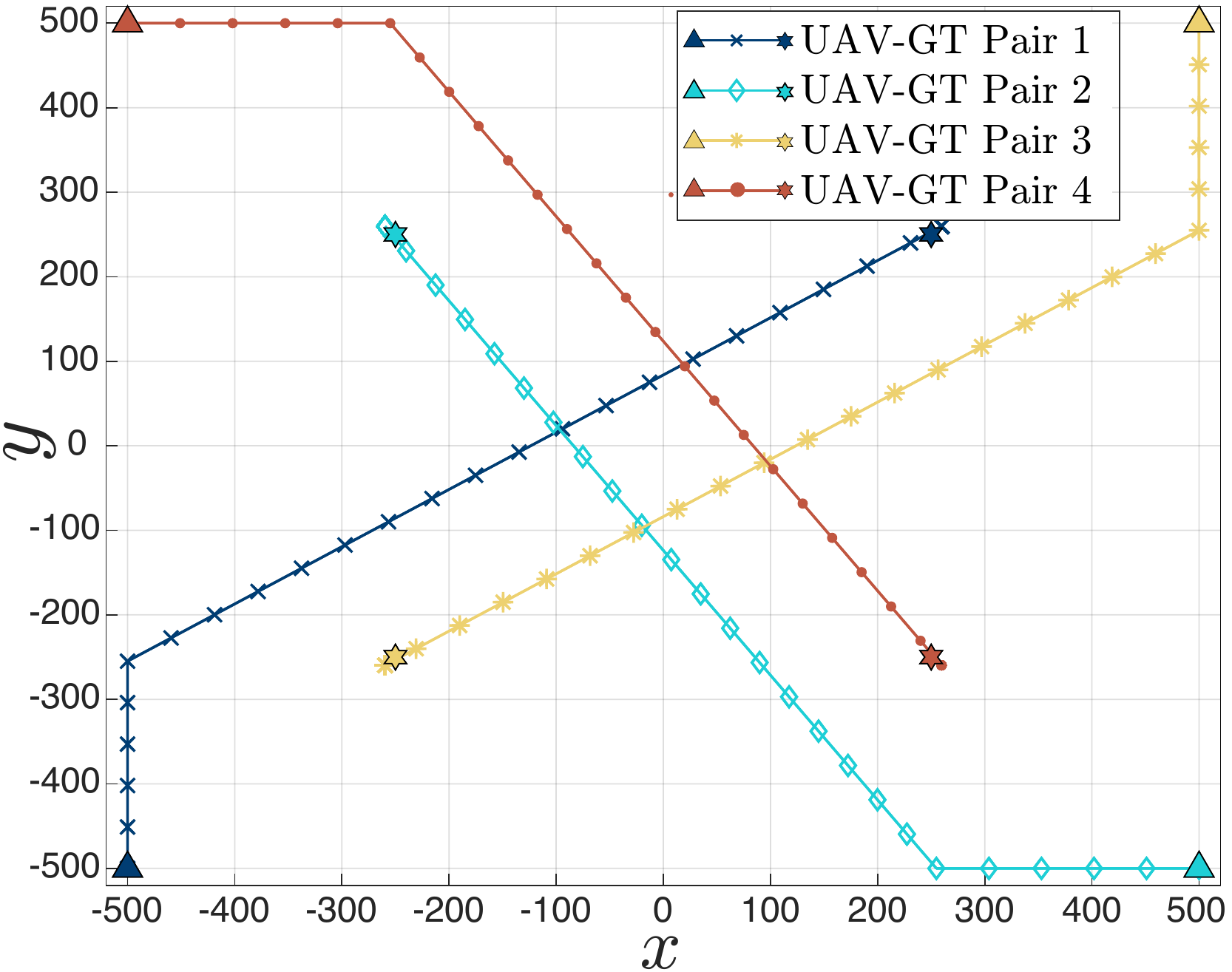}}~
	\subfigure[Optimized trajectories on x-y plane with fixed height.]
			  {\label{fig:OptiTraj:FixedH}\includegraphics[width=0.4\linewidth]{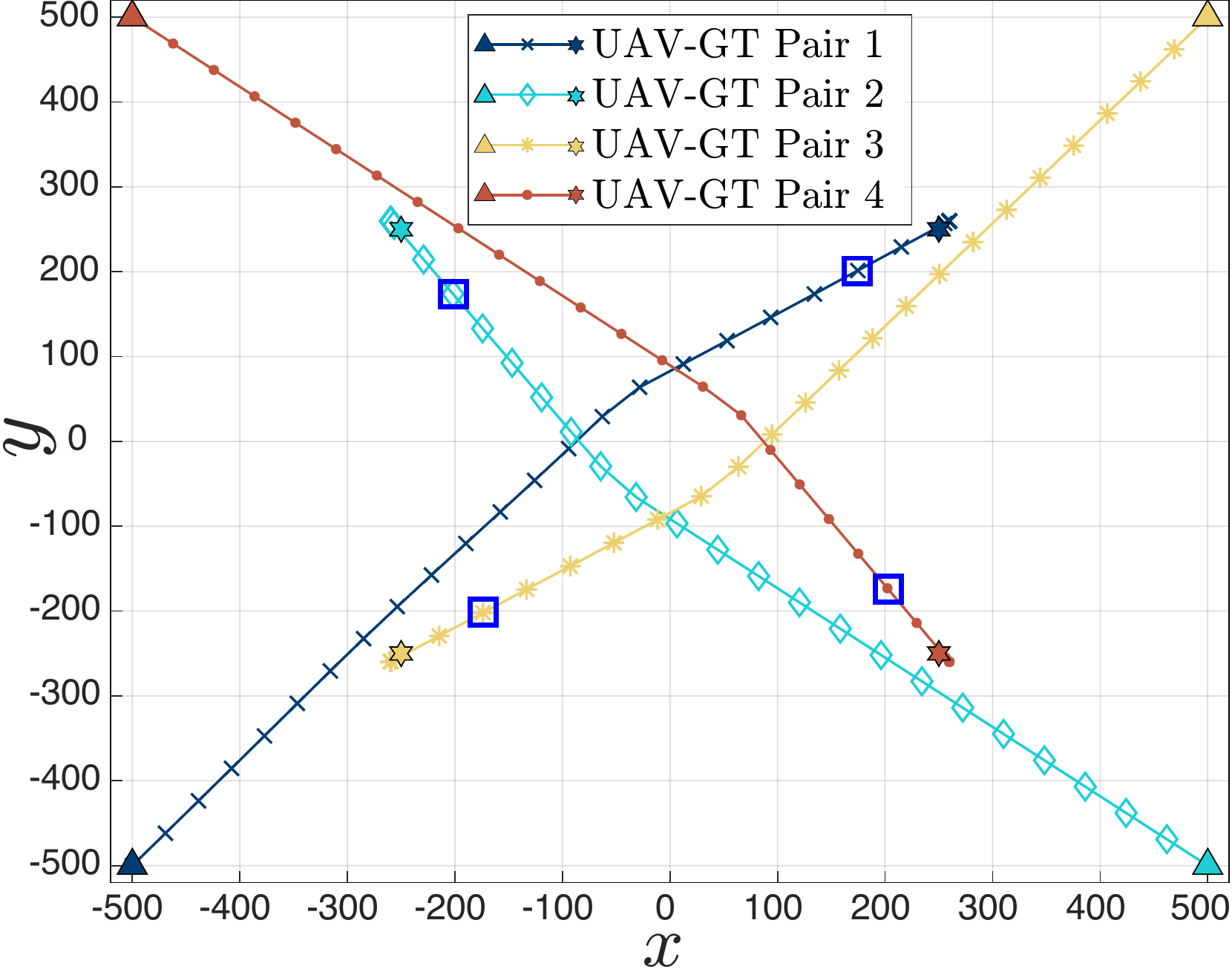}}\\
	\subfigure[Optimized transmission rate with fixed height.]        
			  {\label{fig:OptiRate:FixedH}\includegraphics[width=0.5\linewidth]{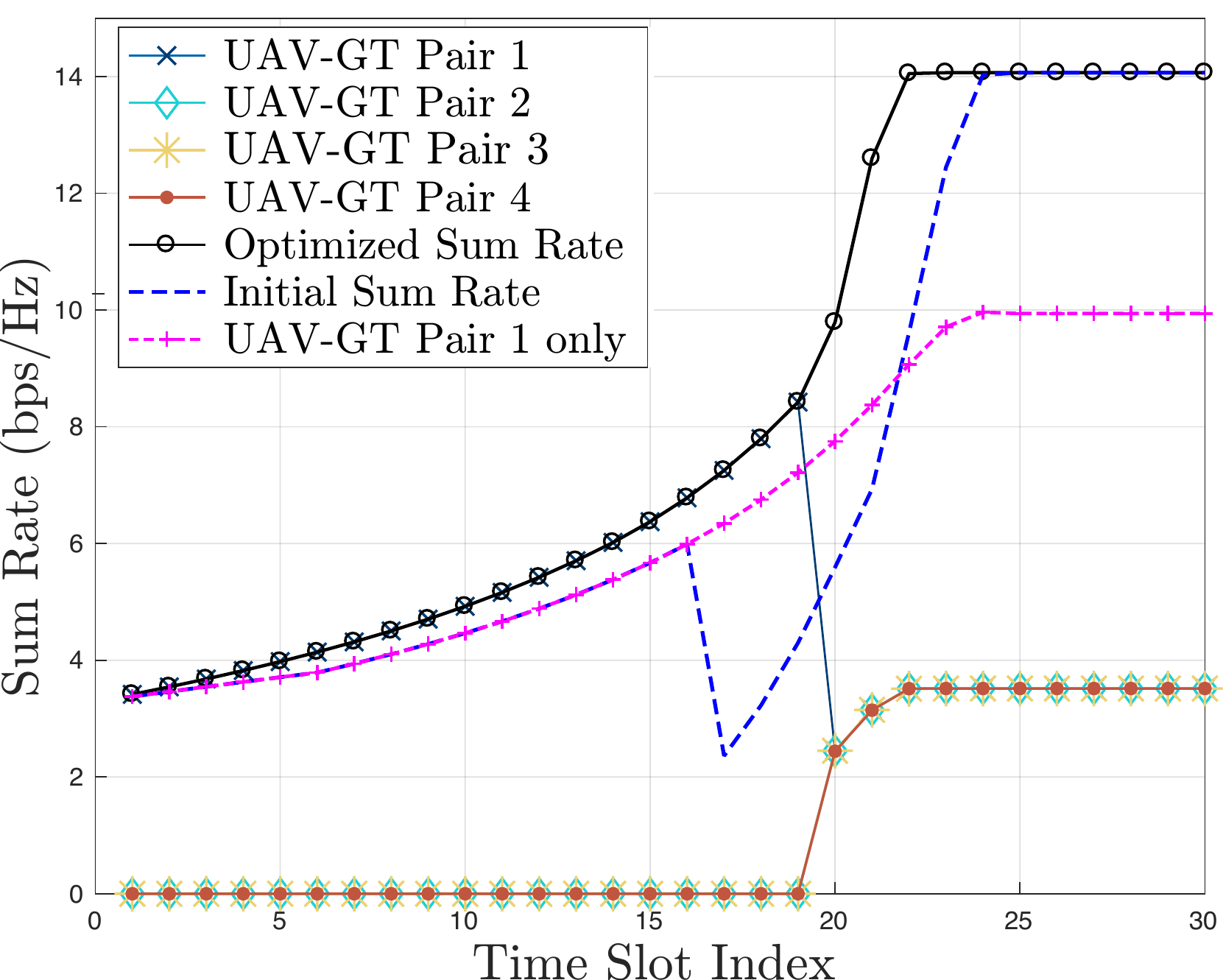}}	\\
	\subfigure[Optimized trajectory on x-y plane.]  
			  {\label{fig::DoF_Height_xy} \includegraphics[width=0.4\linewidth]{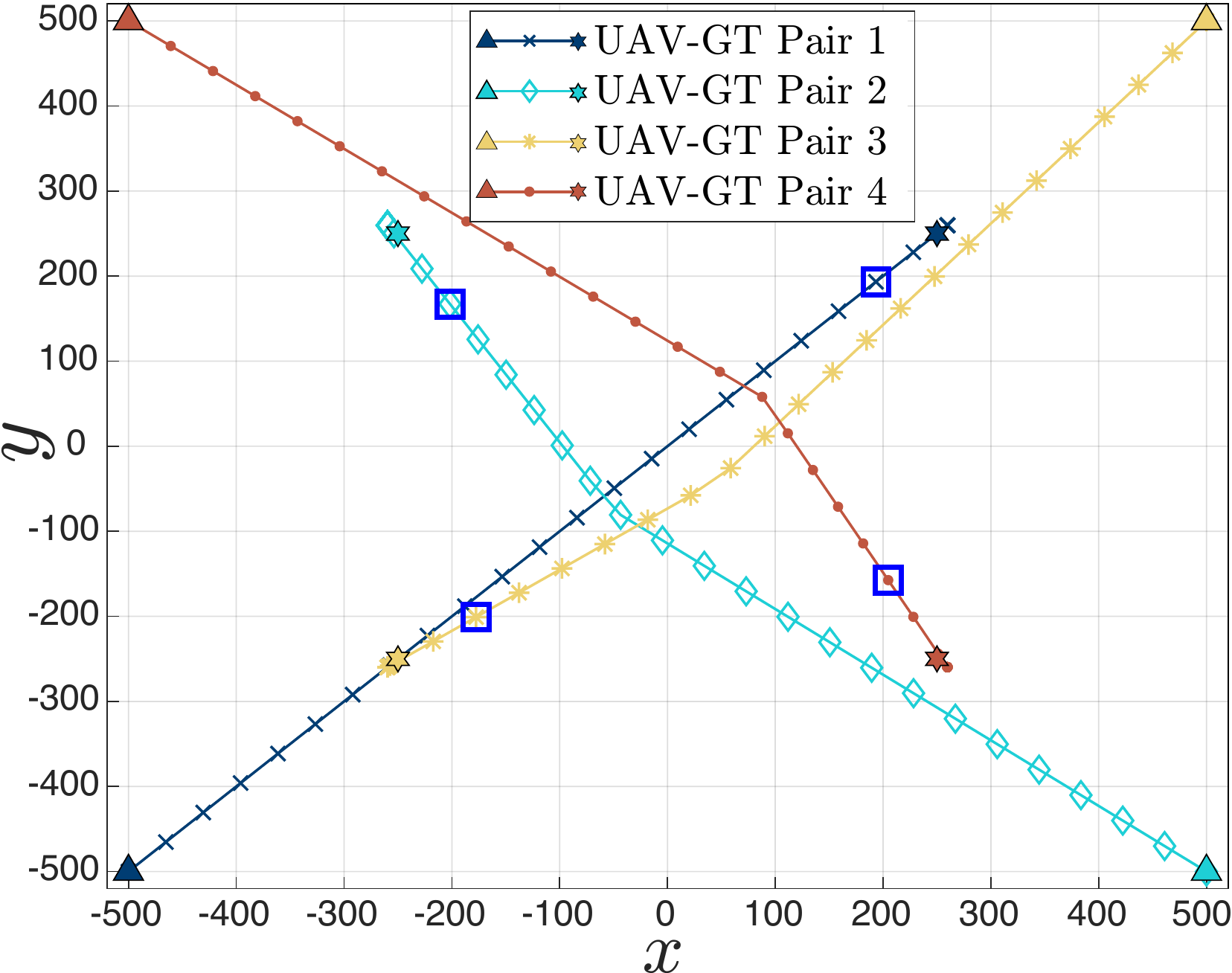}}
	\subfigure[Altitude of the optimized trajectory. ]  
			  {\label{fig::DoF_Height_z} \includegraphics[width=0.41\linewidth]{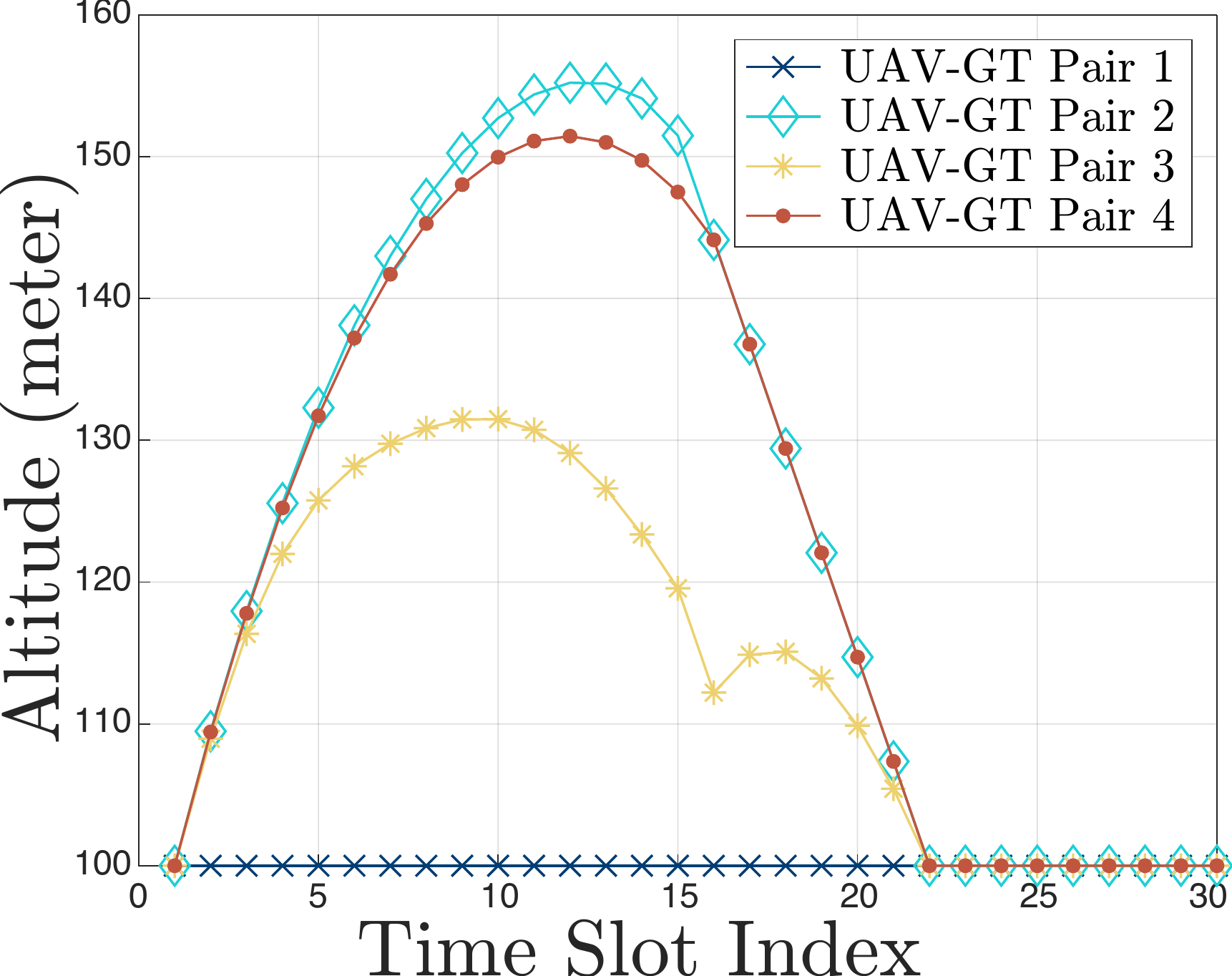}}
	\caption{Trajectories and transmission rates for $K=4$ UAVs, where the flight altitude is fixed at 100 m in (a)-(c), and allowed to vary in (d) and (e), the GTs are initially located at $(\pm250,\pm 250, 0)$, while the associated UAVs are located at $(\mp500,\mp500,100)$, respectively.
	The colored triangles ($\triangle$) and the colored stars ($\bigstar$)  in (a), (b) and (d) denote the initial locations of UAVs and GTs, respectively.
	The blue square markers ({\blue$\square$}) in (b) and (d) denote the UAV locations at the 20th time slot.}
	\label{fig:DoF_Height}
\end{figure}

We first examine how the ability of changing the flying altitude may help the UAVs avoid from collision.
To start with, we consider the case with fixed flying altitude for all UAVs by setting $\Hmax= \Hmin=\dmin=100$ m.
Assume that $K=4$ UAVs are initially located at $(\mp500,\mp500,100)$, respectively, and the GTs are located at $(\pm250,\pm 250, 0)$, respectively.
The hovering locations of the UAVs are first obtained by solving \eqref{Problem::DP}. The initial feasible trajectories of UAVs are shown in Fig. \ref{fig:InitTraj:FixedH} and the trajectories optimized by the proposed Algorithm \ref{Alg2} are shown in Fig. \ref{fig:OptiTraj:FixedH}.
It is observed that the optimized hovering point of each UAV is in close proximity to its own serving GT (not exactly on top of the GT).
Furthermore, due to fixed flying altitude and for collision avoidance, the UAVs have non-trivial flying trajectories.  
Figure  \ref{fig:OptiRate:FixedH} displays the sum rate at different time slots achieved by Algorithm \ref{Alg2} (Optimized Sum Rate) and the sum rate achieved when the UAVs fly with the initial trajectories and use WMMSE-optimized transmission powers (Initial Sum Rate) as well as the rate when there is only UAV 1 present in the network (UAV-GT Pair 1 only). One can clearly observe from the figure that joint TPC optimization can greatly improve the aggregate sum rate performance.
Besides, one can see that only UAV 1 sends information in the first 20 time slots, and then the UAVs start to share the spectrum for transmission.
This implies that the joint TPC optimization allows the UAVs to dynamically switch between TDMA and spectrum sharing modes.

Next, we allow the UAVs to change the flying altitudes and set $\Hmax=500$ m.
The horizontal path and the altitudes of the optimized flying trajectories of the four UAVs are shown in Fig. \ref{fig::DoF_Height_xy} and \ref{fig::DoF_Height_z}, respectively. 
It is observed that different from \ref{fig:OptiRate:FixedH}, the horizontal trajectories of UAVs become more straight, whereas, in order to avoid collision, the UAVs may vary their flying altitudes dynamically.
In particular, UAV 1 establishes a straight-and-level flight to its hovering point whereas other UAVs ascend to different altitudes and gradually descend to approach the respective hovering points. 

The simulation results imply that the extra degree of freedom of flying on different altitudes can greatly help the UAVs for collision avoidance.
Besides, this enables a simple way to find feasible trajectories of UAVs as the initial values for SCA based algorithms, such as the AO method \cite{T:Relay:1S1R1Dzeng,Wu2018TWC} and the proposed Algorithms \ref{Alg1} and \ref{Alg2}, as detailed next.

\subsection{TPC Initialization}\label{Sec::TrajInit}
\begin{figure}[!t]\centering
	\subfigure[Initial trajectory on x-y plane.] 
		{\label{fig:fig1b_init_a}\includegraphics[width=0.4\linewidth]{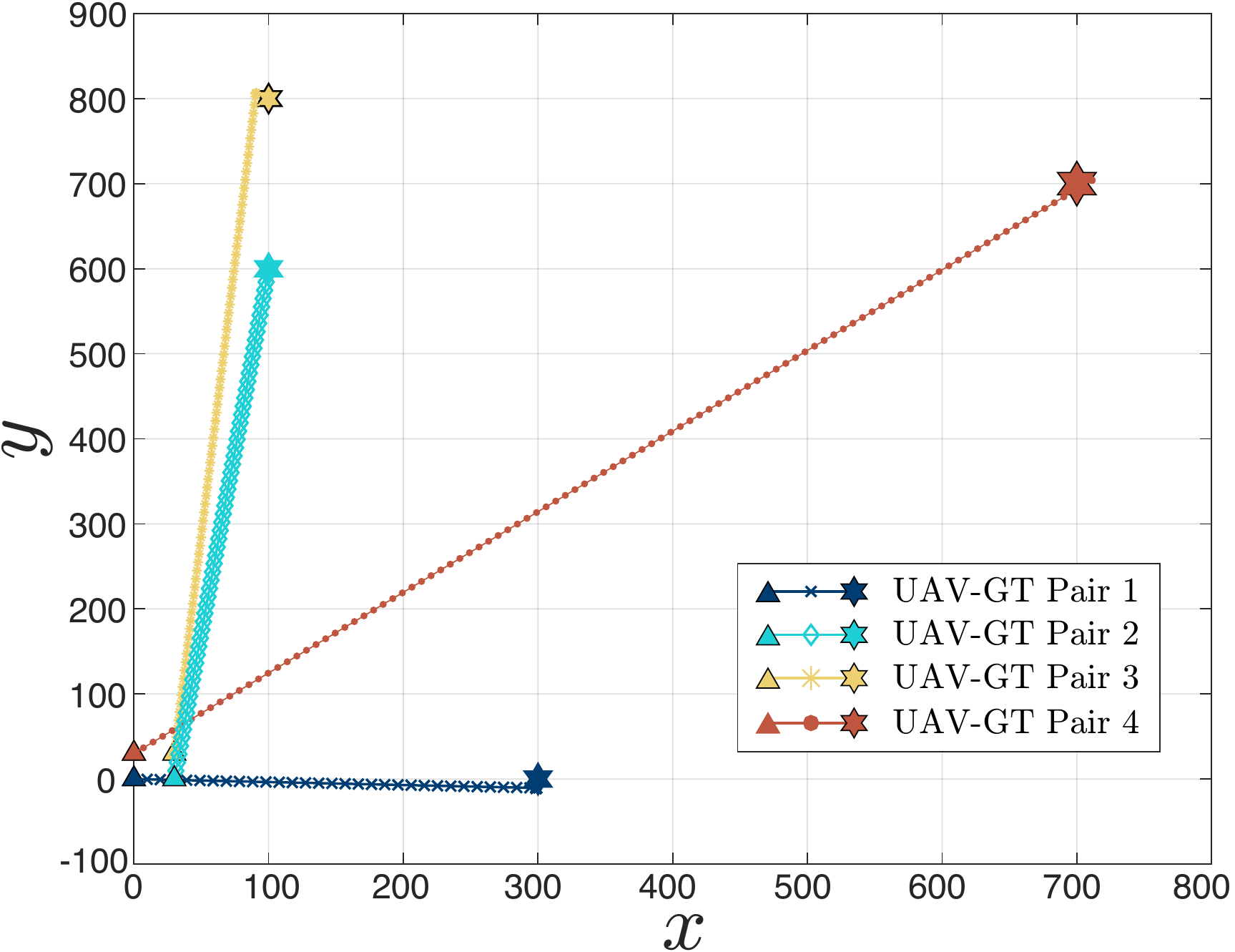}}~~~
	\subfigure[Initial trajectory on flight height vs. timeslot index.]
	{\label{fig:fig1b_init_b} 
		\includegraphics[width=0.4\linewidth]{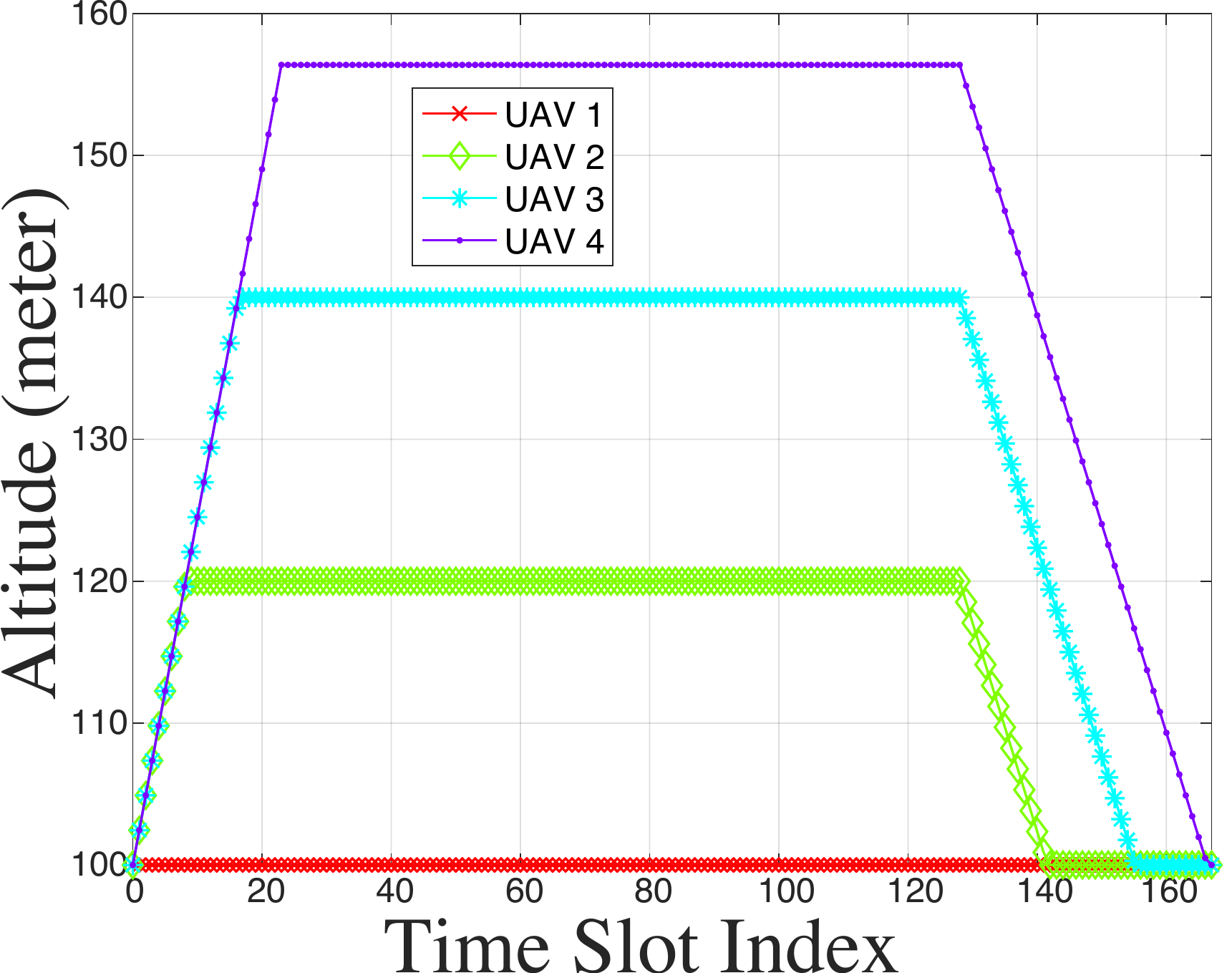}}\\
	\caption{Initial trajectories of UAVs for $K=4$; the initial locations of the UAVs (denoted by $\triangle$ symbols) are $(0,0,100),(30,0,100),(0,30,100),(30,30,100)$, respectively; 
		the GTs are placed at $(300,0,0),(100,600,0),(700,700,0),(100,800,0)$, respectively (denoted by $\bigstar$ symbols).}
	\label{fig:fig1b_init}
\end{figure}

For the deployment optimization problem \eqref{Problem::DP}, the x-y coordinates of the UAV hovering locations are initialized by that of their respective GTs, and the initial altitudes are set to $\Hmin$, namely, $\qb_k^0 = [(\sb_k)_{1:2}^T ~\Hmin]^T$ for all $k$.
For the TPC problems (i.e, \eqref{Problem::SEM3} and FDMA/TDMA formulation \eqref{Problem::SEM:OMA}), a four-step procedure is employed to initialize the UAV trajectories $\{\qb_k[n]\}$. 
In step 1, we obtain a set of hovering positions $\{\qb_k^\star\}$ by solving \eqref{Problem::DP}. 
Since $\qb_k^\star$ for a UAV $k$ with $p_k^\star=0$ is not unique, we reset $(\qb_k^\star)_{1:2}=(\sb_k)_{1:2}$ and $(\qb_k^\star)_3=\Hmin$.
In step 2,  given the initial and hovering locations, each UAV $k$ ascends to the altitude of $\Hmin+(k-1)\dmin$ with the maximum ascending speed and in the meantime flies towards its hovering location $(\qb_k^\star)_{1:2}$ with the maximum level-flight speed, provided that there is no collision with other UAVs.
Then, in step 3, each UAV $k$ flies to its hovering location $(\qb_k^\star)_{1:2}$ by straight and level flight at full speed. 
Finally, in step 4, each UAV $k$ descends to the optimal height $(\qb_k^\star)_3$.
Then the value of $M$ in problem \eqref{Problem::SEM3} can be determined based on the above initial trajectories of UAVs. 
In Fig. \ref{fig:fig1b_init}, we illustrate the initial trajectories of the UAVs for the case of $K=4$.
Given the initial trajectory, the initial transmission powers of UAVs are obtained by the WMMSE algorithm \cite{WMMSE}. 
Note that the initialization procedure described above yields straight trajectories as in Fig. \ref{fig:fig1b_init_a}. 
So the value of $M$ is likely to be smaller, which is desired by Property \ref{Prop:M}.

\subsection{Optimized TPC Solution}
\begin{figure}\centering
	\subfigure[Optimized trajectory on x-y plane.] 		{\label{fig:RT:xy}
		\includegraphics[width=0.4\linewidth]{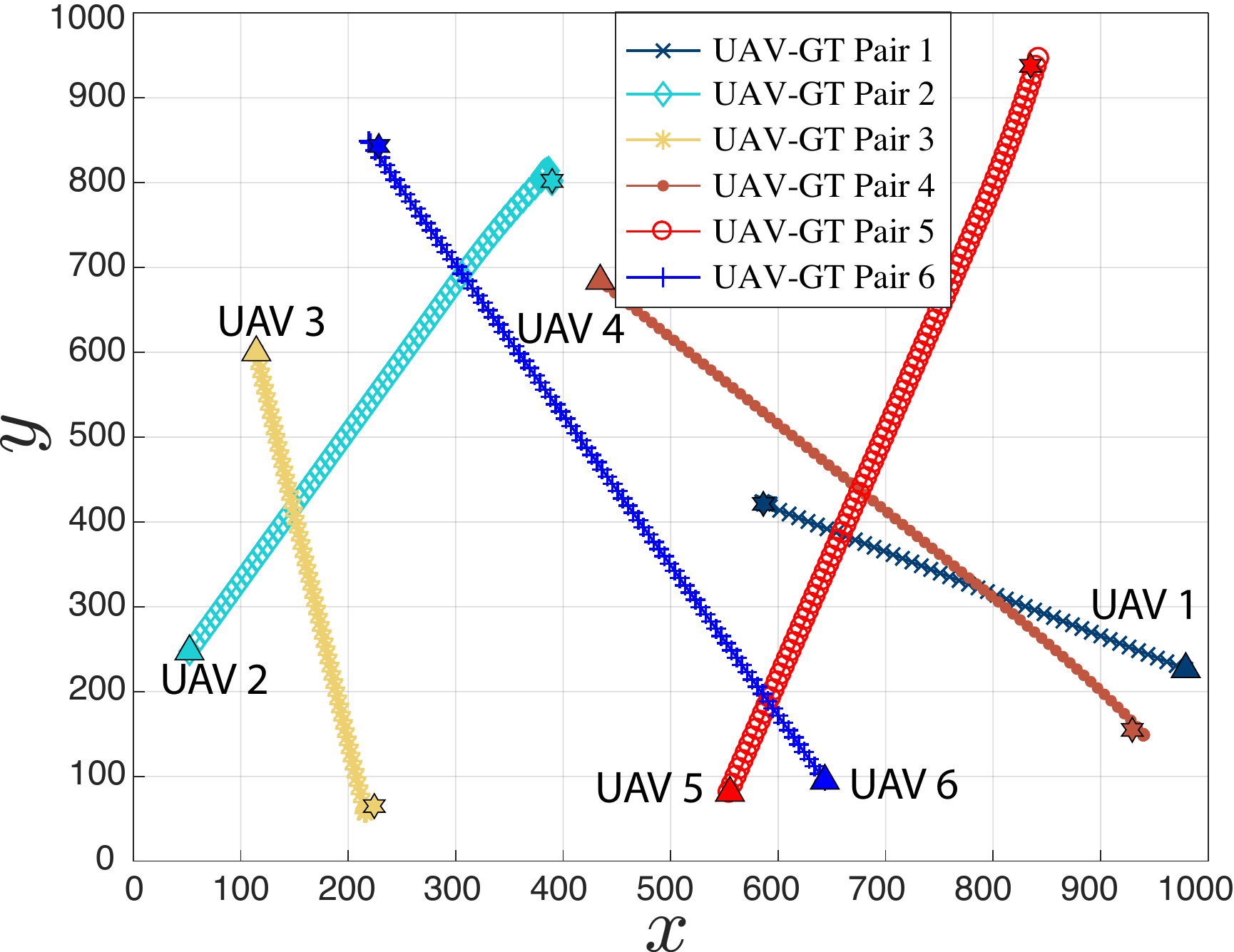}}
	\subfigure[Optimized altitude vs. time slot index.]	{\label{fig:RT:z}
		\includegraphics[width=0.4\linewidth]{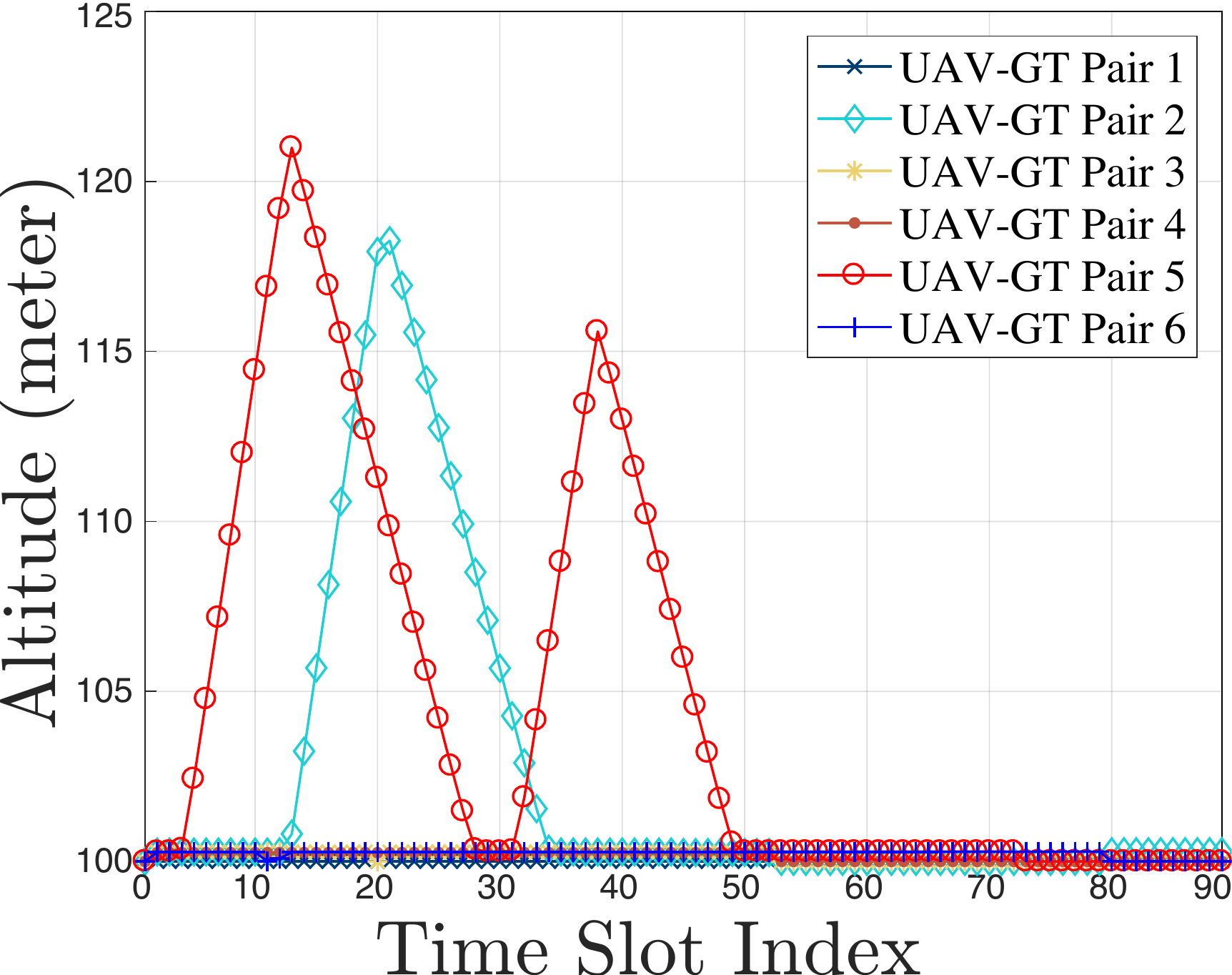}}\\
	\subfigure[Achievable rates vs. time slot index. ]  {\label{fig::RT:R} 	
		\includegraphics[width=0.5\linewidth]{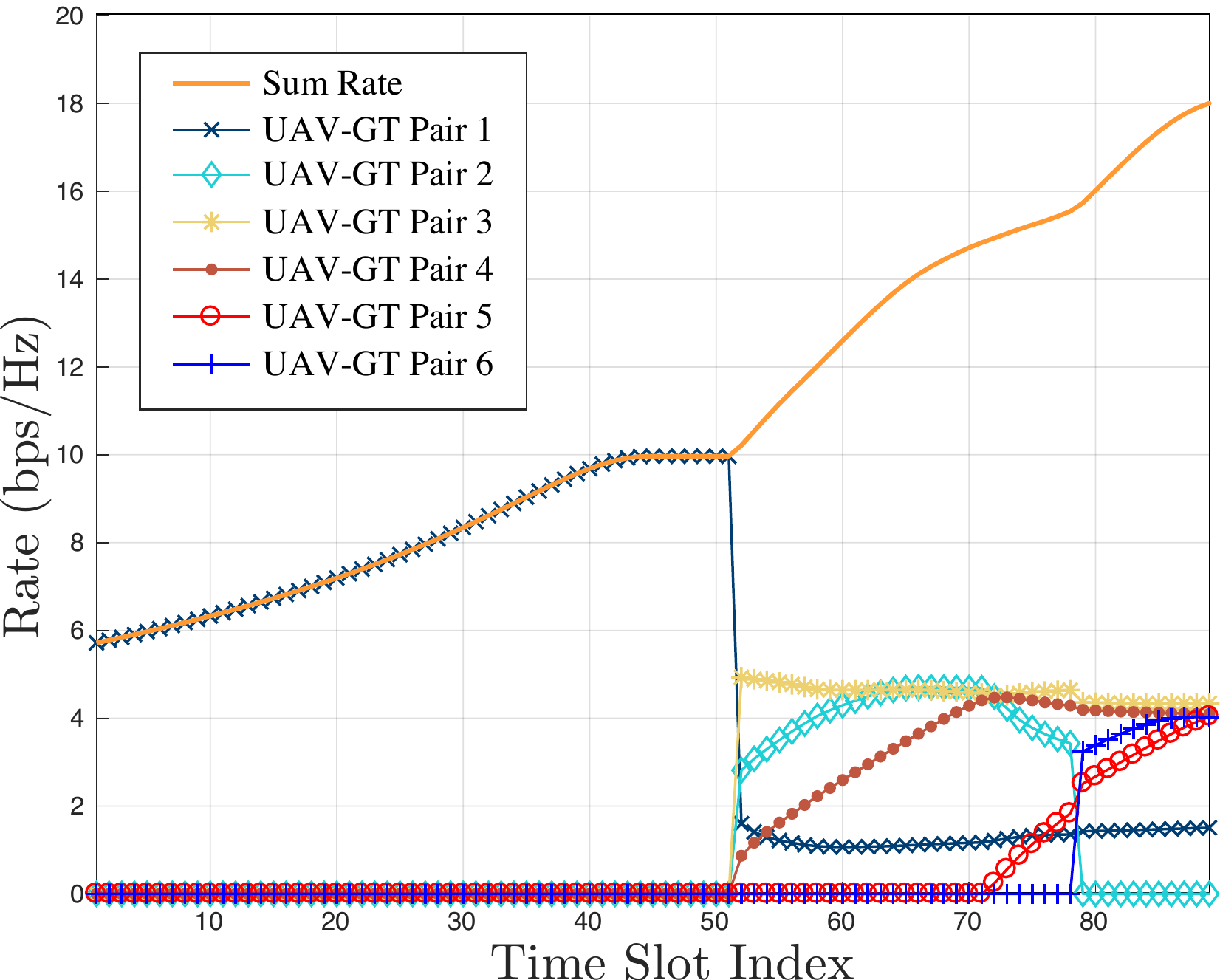}}
	\caption{The optimized trajectory, power and rate for $K=6$ with $\Hmax=500$ m, $\dmax=20$ m and $\Pmax = 30$ dBm.
			 The colored triangles ($\triangle$) and the colored stars ($\bigstar$)  in (a) denote the initial locations of UAVs and GTs, respectively.}
	\label{fig:RandomTopology}
\end{figure}
Fig. \ref{fig:RandomTopology} shows the trajectories, transmission powers and achievable sum rates optimized by the proposed Algorithm \ref{Alg2}
for the case of $K=6$ and $M=92$.  The initial UAVs and GTs positions are  randomly generated within a 1 km $\times$ 1 km square.
One can see from Fig. \ref{fig:RT:xy} that the optimized hovering locations of UAVs are near their serving GTs, and from Fig. \ref{fig:RT:z} that the hovering altitudes are all close to $\Hmin$. 
While all UAVs fly straight to the hovering points with maximum level speed on the x-y-plane, they dynamically change their flying altitudes for collision avoidance.
Specifically, as seen from Fig. \ref{fig:RT:z},  UAV 3 and UAV 5 ascend and then descend so as to avoid collision with UAV 2 and UAV 4, respectively.
It is also observed from Fig. \ref{fig::RT:R} that UAV 2 is chosen to be silent at its hovering location because
the cross-link interference between UAV 2  and UAV 4 as well as between UAV 6 are strong. 
Nonetheless, one can see that it is still possible for UAV 2 to transmit data to its GT during the flight (from time slot 56 to 78) by properly coordinating with nearby UAVs.

\subsection{Convergence of Parallel TPC Algorithm}
\begin{figure}[!t]	\centering
	\subfigure[Realization 1] {\label{fig:convergence_a} \includegraphics[width=0.5\linewidth]{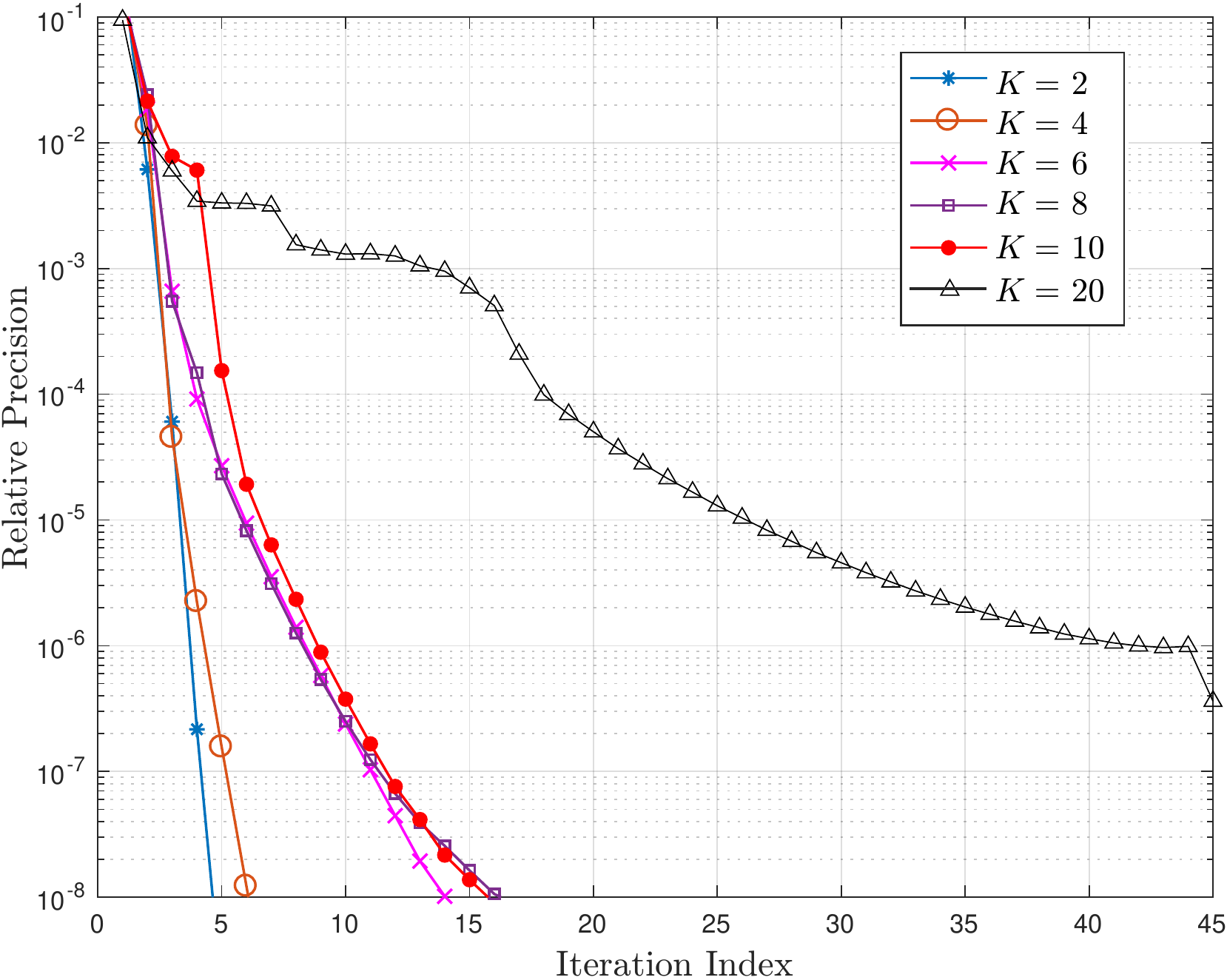} }
	\subfigure[Realization 2] {\label{fig:convergence_b} \includegraphics[width=0.5\linewidth]{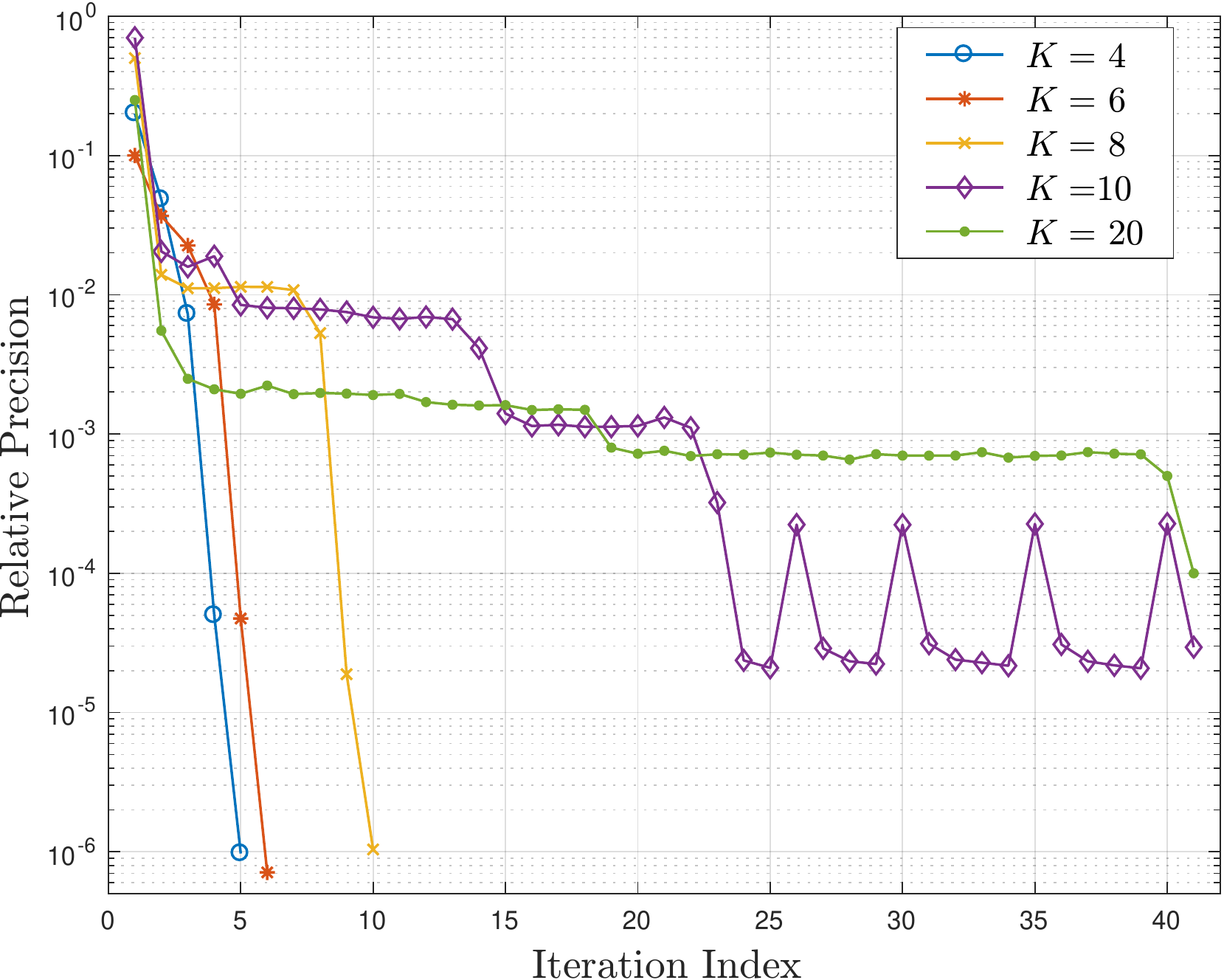} }
	\caption{Typical convergence curves of the proposed parallel TPC algorithm (Algorithm \ref{Alg2}) for $K\in\{2,4,6,8,10,20\}$; the initial locations of UAVs and the locations of GTs are randomly generated. }
	\label{fig:convergence}
\end{figure}
\begin{table*}[!t]\centering
	\caption{Average computation time and aggregate sum rate comparison}
	\label{cputime_random}
	\includegraphics[width=.9\linewidth]{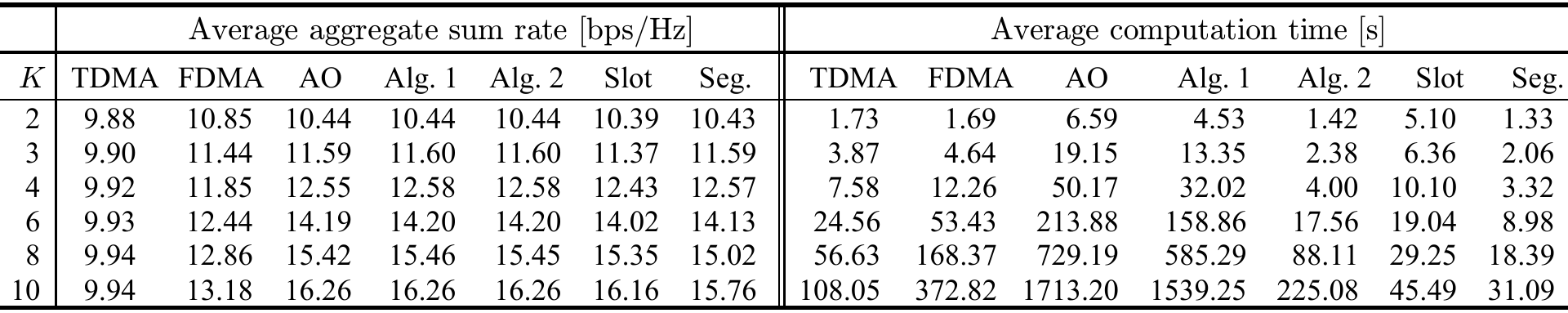}
\end{table*}

Fig. \ref{fig:convergence} presents the typical convergence curves (relative precision\footnote{The relative precision is defined as $\frac{R^r-R^{r-1}}{R^{r-1}}$, where $R^r$ is the aggregate sum rate achieved by $\{a_k^{r+1}[n], \qb_k^{r+1}[n]\}_{n=1}^M$ in the $r$th iteration of Algorithm \ref{Alg2}.} versus the iteration index $r$) of the proposed parallel TPC algorithm (Algorithm \ref{Alg2}), for various numbers of UAVs. 
The initial locations of UAVs and the locations of GTs are randomly generated. Fig. \ref{fig:convergence_a} and Fig. \ref{fig:convergence_b} are obtained by two different random realizations, and they respectively represent two typical convergence behaviors of Algorithm \ref{Alg2}. 
Specifically, one can see from the two figures that Algorithm  \ref{Alg2} can converge to a high-accuracy solution within 10 iterations when the number of UAVs is no greater than 10.
When the number of UAVs increases to 20,  the convergence becomes slower. However, as seen from the two figures, the parallel TPC algorithm (Algorithm \ref{Alg2}) can still achieve a solution with relative precision less than $10^{-3}$ within $30$ iterations.


\subsection{Throughput and Computation Time Comparison}
In this subsection, we examine the computation time and the achieved aggregate sum rate of the proposed algorithms.
Besides the proposed TPC algorithm (Algorithm \ref{Alg1}), parallel TPC algorithm (Algorithm \ref{Alg2}) and the {segment-by-segment} algorithm (Algorithm \ref{Alg3} with $N_{\rm seg}=40$), we also implement the AO method for solving problem \eqref{Problem::SEM3} as well as the FDMA and TDMA schemes presented in Section \ref{sec: FDMA/TDMA}. In addition, a ``slot-by-slot" scheme which corresponds to the segment-by-segment scheme with $N_{\rm seg}=1$ is implemented. The results are presented in Table \ref{cputime_random}. All UAVs are initially located around $(0,0,0)$, and spaced with uniform distance of $d_{\min}=20$ m.
The GTs are initially randomly deployed within a square of length 1 km. 
The results are obtained by averaging over 100 random realizations and each over the flying time $T=10$ minutes, and the simulations were performed on a {desktop} computer with a 4-core 3.40 GHz CPU and 4 GB RAM.

Firstly, one can see that the parallel TPC algorithm (Algorithm \ref{Alg2}) yields almost the same aggregate sum rate as the centralized TPC algorithm (Algorithm \ref{Alg1}), but with significantly reduced computation time, especially when $K\ge 3$. 
It is also observed that the proposed segment-by-segment algorithm (Algorithm \ref{Alg3}) can further reduce the computation time. 
However, as observed from the table, the segment-by-segment algorithm may have about 5\% loss of the aggregate sum rate when compared to Algorithm \ref{Alg1}. The slot-by-slot method is less time efficient than the general segment-by-segment method.
Comparing the proposed TPC algorithms with the AO method, one can see that the AO method achieves almost the same aggregate sum rate performance as the proposed TPC algorithms. 
However, the AO method requires more computation time. 

Lastly, we can observe that
the FDMA scheme can outperform the TDMA scheme for all $K$ and even the non-orthogonal scheme for $K=2$ in terms of the aggregate sum rate.
This is mainly due to the fact that we only consider the short-term peak power constraints instead of the long-term average power limitation for UAV transmission.
However, the FDMA scheme requires more computation time than TDMA for all $K$. Further,  
the non-orthogonal schemes outperform the FDMA scheme, especially when the number of UAV-GT pairs is moderately large ($K \geq 3$). 
This shows the potential advantage of UAV spectrum sharing as long as trajectories and transmission powers of the UAVs are properly coordinated.

\section{Conclusion}\label{sec: conclusion}
In this paper, we have studied the joint TPC design problem for the multiuser UAV-IC. Since the TPC problem is NP-hard and involves a large number of optimization variables,  efficient suboptimal algorithms have been developed in this paper. Specifically, we have shown that for round-trip operation where the UAV needs to return the initial location after the mission, the optimal TPC solutions have a symmetric property and the problem can be approximately solved by first solving the optimal hovering location problem \eqref{Problem::DP} followed by solving the dimension-reduced TPC problem \eqref{Problem::SEM3}.
To find an efficient suboptimal solution of \eqref{Problem::SEM3}, we have proposed an SCA-based algorithm (Algorithm \ref{Alg1}) which, unlike the existing AO based methods, jointly updates the trajectory and transmission power variables in each iteration. For efficient implementation in large scale scenarios with large number of UAVs and/or time slots, we have further proposed the parallel TPC algorithm (Algorithm \ref{Alg2}) and the segment-be-segment method (Algorithm \ref{Alg3}). {Simulation results have shown that the proposed TPC algorithm {achieves almost the same aggregate sum rate performance as the AO method, but requires much less computation time}.} Moreover, the parallel TPC algorithm can achieve nearly the same aggregate sum rate as its centralized counterpart, but with substantially reduced computation time. The segment-be-segment algorithm can further reduce the computation time, though with slight performance degradation. 
The simulation results have also shown that the UAV network with spectrum sharing between different links can outperform the TDMA scheme. 

The current work may motivate several research directions in the future. 
Firstly, as observed from Figure \ref{fig::RT:R}, depending on the relative locations of the UAVs, the UAVs either operate under orthogonal time slots or share the spectrum. It is therefore interesting to consider a general scheme that allows the UAVs to dynamically switch between FDMA and spectrum sharing under both peak and average power constraints. 
Secondly, it is important to extend the current work to scenarios with multi-antenna GTs and with numbers of GTs larger than the UAVs \cite{Liuliang_BF_UAV}, i.e., the broadcast or multicast channels. 

\footnotesize
\bibliographystyle{IEEEtran}
\bibliography{UAVref}

\end{document}